\documentclass[12pt]{amsart}
\usepackage{fullpage}
\usepackage{color}
\usepackage{pstricks,pst-node,pst-plot} 
\usepackage{graphicx,psfrag} 
\usepackage{color} 
\usepackage{tikz}
\usetikzlibrary{decorations.markings, hobby, knots}
\usepackage{pgffor}
\usepackage{hyperref}
\usepackage{todonotes}
\usepackage{subfigure}
\usepackage{mathtools}
\usepackage{verbatim}
\usepackage[utf8x]{inputenc} 
\usepackage{enumerate}

\usepackage{bm}
\usepackage{multirow}

\newtheorem{theorem}{Theorem}[section]
\newtheorem*{theorem-non}{Theorem}
\newtheorem{lemma}[theorem]{Lemma}

\newtheorem{proposition}[theorem]{Proposition}

\theoremstyle{definition}
\newtheorem{definition}[theorem]{Definition}
\newtheorem{example}[theorem]{Example}

\newtheorem{notation}[theorem]{Notation}

\theoremstyle{remark}
\newtheorem{remark}[theorem]{Remark}

\numberwithin{equation}{section}

\definecolor{gray}{rgb}{.5,.5,.5}

\definecolor{black}{rgb}{0,0,0}

\definecolor{blue}{rgb}{0,0,1}

\definecolor{red}{rgb}{1,0,0}

\definecolor{green}{rgb}{0,1,0}

\definecolor{yellow}{rgb}{1,1,.4}

\newrgbcolor{purple}{.5 0 .5}
\newrgbcolor{black}{0 0 0}
\newrgbcolor{white}{1 1 1}
\newrgbcolor{gold}{.5 .5 .2}
\newrgbcolor{darkgreen}{0 .5 0}
\newrgbcolor{gray}{.5 .5 .5}
\newrgbcolor{lightgray}{.75 .75 .75}
\newrgbcolor{lightred}{.75 0 0}

\newcommand{\crossneg}{
\begin{tikzpicture}[baseline=-2]
\draw[white,line width=1.5pt,double=black,double distance=.5pt] (0,-0.1) -- (0.3,0.2);
\draw[white,line width=1.5pt,double=black,double distance=.5pt] (0,0.2) -- (0.3,-0.1);
\end{tikzpicture}}

\newcommand{\crosspos}{
\begin{tikzpicture}[baseline=-2]
\draw[white,line width=1.5pt,double=black,double distance=.5pt] (0,0.2) -- (0.3,-0.1);
\draw[white,line width=1.5pt,double=black,double distance=.5pt] (0,-0.1) -- (0.3,0.2);
\end{tikzpicture}}

\newcommand{\crossposor}{
\begin{tikzpicture}[baseline=-2]
\draw[white,line width=1.5pt,double=black,double distance=.5pt] (0,0.2) -- (0.3,-0.1);
\draw[<-] (0,0.2) -- (0.3,-0.1);
\draw[white,line width=1.5pt,double=black,double distance=.5pt] (0,-0.1) -- (0.3,0.2);
\draw[->] (0,-0.1) -- (0.3,0.2);
\end{tikzpicture}}

\newcommand{\crossnegor}{
\begin{tikzpicture}[baseline=-2]
\draw[white,line width=1.5pt,double=black,double distance=.5pt] (0,-0.1) -- (0.3,0.2);
\draw[->] (0,-0.1) -- (0.3,0.2);
\draw[white,line width=1.5pt,double=black,double distance=.5pt] (0,0.2) -- (0.3,-0.1);
\draw[<-] (0,0.2) -- (0.3,-0.1);
\end{tikzpicture}}

\newcommand{\crosssingor}{
\begin{tikzpicture}[baseline=-2]
\draw[<-] (0,0.2) -- (0.3,-0.1);
\draw[->] (0,-0.1) -- (0.3,0.2);
\end{tikzpicture}}

\begin{document}

\title{Average Casson invariant of 2-bridge knots}

\author{Moshe Cohen}
\address{Department of Mathematics, State University of New York at New Paltz, New Paltz, NY 12561}
\email{cohenm@newpaltz.edu}

\author{Marina Ville}
\address{Univ. Paris Est Creteil, CNRS, LAMA, F-94010 Creteil, France}
\email{villemarina@yahoo.fr}

\begin{abstract}
The Casson invariant of a knot is a Vassiliev or finite type invariant of type 2.  We show that the average or expected absolute value of the Casson invariants of all 2-bridge knots with crossing number $c$ approaches $c/6+1/36 - (9/4)(-1)^c$ as $c$ approaches infinity.  We compare this result to some in the random knot literature.  We conclude with two infinite families of knots, the $T(2,c)$ torus knots and a subclass of the double twist knots,  whose Casson invariants are maximal and minimal, respectively.
\end{abstract}

\subjclass[2020]{57K10}
\keywords{rational knot}

\maketitle

\section{Introduction}

In recent years, the first author with Lowrance and others and, independently, Suzuki and Tran,  have published a series of similar papers asking the question:  What does the expected value of a knot invariant look like when a 2-bridge knot of a given crossing number $c$ is selected at random?  Invariants considered so far include Seifert (3-)genus \cite{Coh:lower, CohLow, CohLowURSI, SuzukiTran, RayDiao}, braid index \cite{LowURSI:braids, SuzTra:braids}, signature and 4-genus \cite{CohLowMadRaa, CohMadAYURE} with Madras and others, and unoriented genus and crosscap number \cite{CohKinLowShaVan} with Kindred, Shanahan, and Van Cott.  
The main technique in some of these works involves a well-behaved recursion on the partially doubly-counted set $T(c)$ of 2-bridge knots with crossing number $c$.  Together with the set $T_p(c)$ of palindromic 2-bridge knots with crossing number $c$, these count exactly twice every knot in the set $\mathcal{K}(c)$ of 2-bridge knots with crossing number $c$.

In this work we address Vassiliev or finite-type invariants first introduced by Vassiliev \cite{Vas} and also studied by Kontsevich \cite{Kont:Vas}, Polyak and Viro \cite{PolVir:Vas} using Gauss diagrams, Bar-Natan \cite{BN:Vas}, and Polyak \cite{Pol:inv}.  These satisfy the Vassiliev skein relation $v(\crosssingor)=v(\crossposor)-v(\crossnegor)$ on singular knots and links.  Such an invariant is of type $k$ if it vanishes for links with more than $k$ singular points.  For more details see the textbook on the subject by Chmutov, Duzhin, and Mostovoy \cite{ChmDuzMos}.

The simplest Vassiliev invariant is the linking number.  The next is $v_2$ which has type 2, coincides with the coefficient of $z^2$ in the Conway polynomial, and reduces modulo 2 to the Arf invariant.  It is also called the Casson invariant of knots
because it appears in the theory of the Casson invariant of homology 3-spheres (as in Morita \cite{Mor:Cas}, Saveliev \cite{Sav:Casson}, and Polyak and Viro \cite{PolVir:Cas}).

We aim to understand how these invariants grow as our knots get larger.  Stoimenow \cite{Sto:thesis} studied asymptotics of these Vassiliev invariants, Willerton \cite{Wil:pairs} studied pairs ($v_2,v_3$) of Vassiliev invariants, Koseleff and Pecker \cite{KosPec:Alex} bounded the Casson invariant from above by the coefficient of the largest power in the Conway polynomial, Abe \cite{Abe:Vas} studied other bounds, and Taniyama \cite{Tan:pairs} studied pairs of invariants (including these) more generally.

Kalfagianni and Lin \cite{KalLin} used Vassiliev invariants as obstructions to finding regular Seifert surfaces with simple complements.  Recently Klotz and Estabrooks \cite{KloEst} used the Casson invariant to distinguish knotted polymers, and Bar-Natan, Bar-Natan, Halacheva, and Scherich \cite{BNBNHS} studied efficient computation of these invariants.

Due to their combinatorial nature, Vassiliev invariants can be easier to compute than other knot invariant, and the Casson invariant $v_2$ is a Vassiliev invariant of type 2 so it is easier to compute than those with higher type.  As such, this one has appeared in recent random models for knot diagrams by Even-Zohar, Hass, Linial, and Nowik \cite{EZHLN} and by the second author \cite{Ville:Casson}.  See Section \ref{sec:random}.

\subsection*{Main result}

We define the sums of Casson invariants over the partially doubly-counted set $T(c)$ of 2-bridge knots with crossing number $c$ and the set $T_p(c)$ of palindromic 2-bridge knots with crossing number $c$ mentioned at the start of the paper.

    The \textbf{total Casson invariant} $v(c)$ from $T(c)$ is defined as 
    \[v(c) = \sum_{w\in T(c)} v_2(K_w),\]
    and the \textbf{total palindromic Casson invariant} $v_p(c)$ from $T_p(c)$ is defined as 
    \[v_p(c) = \sum_{w\in T_p(c)} v_2(K_w).\]

With this notation, the average Casson invariant $\overline{v_2}(c)$ can be computed as
\[\overline{v_2}(c) = \frac{v(c) + v_p(c)}{2|\mathcal{K}(c)|}.\]
Within this context, we present the main result of this paper:

\begin{theorem}
\label{thm:main}
The average Casson invariant of a 2-bridge knot of crossing number $c$ 
is
$$\overline{v_2}(c)\approx\dfrac{v(c)}{2|\mathcal{K}(c)|}\approx\dfrac{(c2^{c-3})/9+2^{c-3}/54-3\cdot2^{c-3}(-1)^c/2}{2(2^{c-3})/3}\approx\dfrac{c}{6}+\dfrac{1}{36}-\dfrac{9}{4}(-1)^c,$$
where we use $\approx$ to denote equivalence as $c$ approaches infinity.
\end{theorem}

We provide the brief proof of this result now as a way to highlight the organization of the rest of the paper.

\begin{proof}
We first rely on Proposition \ref{prop:ignorepalindromes} (which imitates a proof in \cite[Section 4]{CohLowMadRaa}) to show that we can ignore the contribution $v_p(c)$ to the total Casson invariant from the palindromes in $T_p(c)$ and only consider the contribution $v(c)$ to the total Casson invariant from the partially doubly-counted set $T(c)$.

Then the numerator is obtained from Theorem \ref{thm:v}, which states that 
$$v(c) \;  = \; \frac{1}{72}c2^{c}  +\frac{1}{432}2^{c}   -\frac{3}{16}(-2)^{c} + \frac{2}{9}(-1)^cc   +\frac{5}{27}(-1)^{c},$$
and the denominator follows from Theorem \ref{thm:ernstsumners} by Ernst and Sumners counting the number $\mathcal{K}(c)$ of 2-bridge knots with crossing number $c$.
\end{proof}

The formula for $v(c)$ above coincides with data obtained from online databases:  KnotInfo \cite{knotinfo} to determine which knots appear in $T_p(c)$ (and thus only count once in $T(c)$) and the Knot Atlas \cite{knotatlas} to obtain the Casson invariant.  See Table \ref{tab:vexperimental}.  The final row gives the first approximation from Theorem \ref{thm:main}.  Computations performed by the authors indicate that by $c=98$, the difference between the first approximation and the last approximation given in the Main Theorem \ref{thm:main} is on the order of $10^{-14}$.

\begin{table}[h!]
\begin{tabular}{|c|c|c|c|c|c|c|c|c|}
\hline
$c$ &   3 &     4 &     5  &    6 &     7 &     8 &     9   & 10    \\
\hline
$v(c)$ &   1 &     -1 &     7  &    -5 &     35 &     -17 &     159   & -45    \\
\hline
$2|\mathcal{K}(c)|$ &   2 &     2 &     4  &    6 &     14 &     24 &     48   & 90    \\
\hline
$\displaystyle\frac{v(c)}{2|\mathcal{K}(c)|}$ &   $\frac{1}{2}$ &     -$\frac{1}{2}$ &     $\frac{7}{4}$  &    -$\frac{5}{6}$ &     $\frac{5}{2}$ &     -$\frac{17}{24}$ &     $\frac{53}{16}$   & -$\frac{1}{2}$    \\
\hline
\end{tabular}
    \caption{\label{tab:vexperimental} The data in the $v(c)$ row  
    from KnotInfo \cite{knotinfo}  
    and Knot Atlas \cite{knotatlas}  
    that coincide with the formula above, along with Ernst and Sumner's count of 2-bridge knots and our first approximation from Theorem \ref{thm:main}}
\end{table}

\subsection*{Organization}  
We use Section \ref{sec:random} to compare our Main Theorem \ref{thm:main} with results from the literature on random knots.

We provide enough background in Section \ref{sec:background} to understand the notation and ingredients in the Main Theorem \ref{thm:main} above.

The goal of Section \ref{sec:main} is to yield Theorem \ref{thm:v} on the contribution to the total Casson invariant from knots in the partially doubly-counted set $T(c)$.

We build up to Proposition \ref{prop:ignorepalindromes} on the insignificant contribution from palindromes in Section \ref{sec:palindromes} and then conclude with some discussion about maximal and minimal cases since this idea is used in the proof of the proposition.  In particular, we refer to, respectively, $T(2,c)$ torus knots and a subclass of double twist knots that have also been recently studied by, for example, Hoste, Shanahan, and Van Cott \cite{HosShaVanCot:double} and Singh et al 
\cite{SCDD:double}.

\subsection*{Acknowledgments}  We are grateful to the Banff International Research Station for hosting the authors at the Knot Theory Informed by Random Models and Experimental Data (24w5217) workshop and to the organizers of the workshop.  This project was conceived in the beautiful dining room after a meal.  The first author would also like to thank Pierre Berger and Pierre-Vincent Koseleff for the invitation to visit the Institut de Math\'ematiques de Jussieu in Paris as well as Vera Toni\'c for the invitation to visit the University of Rijeka in Croatia, during which portions of this work were completed.

\subsection*{AI acknowledgment}  This research was completed without the use of AI tools.

\section{Comparison of results with literature on random knots}
\label{sec:random}

One motivation for undertaking this work has been to obtain a baseline for invariants of knots that can be used to compare with results obtained from random models.

In Subsection \ref{subsec:randpetal} we compare with random Petaluma knot diagrams.  In Subsection \ref{subsec:randrosette} we compare with random rosette and quasitoric braids.  In Subsection \ref{subsec:randlissajous} we compare with random Lissajous knot diagrams.  In Subsection \ref{subsec:rand3braids} we compare with random 3-braids.

\subsection{Petaluma knot diagrams}
\label{subsec:randpetal}
Even-Zohar, Hass, Linial, and Nowik \cite{EZHLN,EZHLN:distpetal} explored a model of random Petaluma knot diagrams that were developed by Adams et al \cite{Adams:petal, Adams:petal2, Adams:petal3} and also studied by Colton, Glover, Hughes, and Sandberg \cite{CGHS}. They  consider diagrams that have $2n+1$ petals passing through a single ``\"ubercrossing" and a permutation $\pi\in S_{2n+1}$ drawn uniformly at random.
\begin{theorem}[\cite{EZHLN} Theorem 2 restricted to $k=1$]
    \label{thm:Petaluma} 
    The expected value of the Casson invariant of a randomly chosen petal diagram with $2n+1$ petals is a polynomial in $n$ of degree $2$. In particular 
    $$\displaystyle\lim_{n\rightarrow \infty} \dfrac{E(v_2)}{n^2}=\dfrac{1}{24}.$$
\end{theorem}

  The ``\"ubercrossing" accounts for $\binom{2n+1}{2}=2n^2+n$ usual crossings in a usual knot diagram, giving an upper bound for the crossing number of the knot type $c(K)\leq 2n^2+n$. However, since the distribution of knots appearing in this model is still unknown, this result only provides a lower bound for the expected Casson invariant with respect to crossing number of the knot.  
For the crossing number $c$ of such a knot we have
\begin{equation*}
    \lim_{c\rightarrow\infty} \frac{E(v_2)}{c} \geq \lim_{n\rightarrow\infty} \frac{E(v_2)}{2n^2+n} = \frac{1}{48}.
\end{equation*} 

Comparing our Main Theorem \ref{thm:main} 
to Theorem \ref{thm:Petaluma} by \cite{EZHLN}, we do indeed asymptotically get a polynomial that is linear with respect to the crossing number.  Perhaps this comparison can be used to inform on the distribution of crossing numbers of knots in the random petal model, an open and difficult question.  See also Vernitski \cite{Vern}.

\subsection{Rosette and quasitoric braids}
\label{subsec:randrosette}
The second author \cite{Ville:Casson} explored several other models of random knot diagrams.  

The first are closures of rosette $s$-braids $\mathcal{R}(s,k)$ of type $(s,k)$ defined by Lamm \cite{Lam:fourier}.  These braids 
are $k$ factors of the form $\left[ \prod_{\text{odd }j<s}\sigma_j^{\pm1} \prod_{\text{even }j<s}\sigma_j^{\pm1} \right]$.  The second are closures of quasitoric $s$-braids $\mathcal{Q}(s,k)$ of type $(s,k)$ defined by Manturov \cite{Man:quasi}.  These braids 
are $k$ factors of the form $[\sigma_1^{\pm1}\ldots\sigma_{s-1}^{\pm1}]$. Both have $\ell=k(s-1)$ crossings with $c(K)\leq k(s-1)$.  

Both of these diagrams allow for a random model by choosing each crossing to be either $\sigma_i$ or $\sigma_i^{-1}$ with equal probability, a technique also used to study random 2-bridge knots by the first author with Krishnan and Even-Zohar \cite{CoKr,CoEZKr}.

\begin{theorem}
\label{thm:rosette}
    \cite[Theorem 2]{Ville:Casson}  The expected value of the Casson invariant of the closure of a random rosette or quasitoric braid of type $(s,ns+1)$ with braid length $\ell=(ns+1)(s-1)$ crossings and braid index $s$ is $\frac{s(s^2-1)n}{24}$.  
    In particular
    \begin{equation*}
        \lim_{n\rightarrow\infty} \frac{E(v_2)}{n} = \frac{s+1}{24}.
    \end{equation*}
The expected value of the Casson invariant of the closure of a random rosette or quasitoric braid of type $(s,ns-1)$ with braid length $\ell=(ns-1)(s-1)$ crossings and braid index $s$ is $\frac{s(s-1)(5-s)n}{24}$.  
    In particular
    \begin{equation*}
        \lim_{n\rightarrow\infty} \frac{E(v_2)}{n} = \frac{5-s}{24}.
    \end{equation*}
\end{theorem}

Our 2-bridge knots can be written as plat closures of braids on 3 strands.  In the result above for the braid closures for $s=3$ we obtain $\frac{1}{6}$ in the first limit and $\frac{1}{12}$ in the second.  Comparing our Main Theorem \ref{thm:main} to Theorem \ref{thm:rosette} by \cite{Ville:Casson}, we remarkably obtain the same linear coefficient for the first type!  This is interesting despite, again, the crossing number of such a knot being bounded above by $\ell$.  See also the result of Clark, Frank, and Lowrance \cite[Theorem 1.2]{LowURSI:braids} on the average braid index being $\lceil \frac{c}{3} \rceil$ for a 2-bridge knot with crossing number $c$.

\subsection{Lissajous knots}
\label{subsec:randlissajous}
The second author \cite{Ville:Casson} also considered $(p,q)$ checkerboard diagrams, which are plat closures of rosette braids of type $(2q,(p-1))$ with an extra half-factor $\left[ \prod_{\text{odd }j<q}\sigma_j^{\pm1}\right]$ appended at the end.  These are Lissajous knots or cubic billiard knots studied by Jones and Pryzytycki \cite{JonPrz} and  Bogle, Hearst, Jones, and Stoilov \cite{BHJS} and generalized by Koseleff and Pecker \cite{KosPec:billiard}.  These have $2pq-p-q$ crossings with $c(K)\leq 2pq-p-q$.  

These diagrams again allow for a random model by choosing each crossing to be either $\sigma_i$ or $\sigma_i^{-1}$ with equal probability.

\begin{theorem}
\label{thm:Lissajous}
    \cite[Theorem 3]{Ville:Casson}  The expected value of the Casson invariant of a random $(p,q)$ Lissajous knot with braid length $\ell=2pq-p-q$ crossings and braid index $2q$ is $\frac{(p-1)(q-1)}{8}$.  
    In particular
    \begin{equation*}
        \lim_{p\rightarrow\infty} \frac{E(v_2)}{p} = \frac{q-1}{8}.
    \end{equation*}
\end{theorem}

Our 2-bridge knots can be written as plat closures of braids on 3 strands or 4 strands; see work by Koseleff, Pecker, and Rouillier \cite{KosPecRou:10crossings}.  In the result above for the braid closures for $q=2$ we obtain $\frac{1}{8}$ in the limit.

\subsection{3-braids}
\label{subsec:rand3braids}
The second author \cite{Ville:Casson} also studied random 3-braids via a random walk on the braid group of fixed length, having previously used this model with Soret \cite{SorVil}.  It also appeared in work by Dunfield and Tiozzo \cite{DunTio}.
    
    \begin{theorem}
\label{thm:3braid}
\cite[Proposition 3]{Ville:Casson}
The expected value of the Casson invariant of the closure of a random 3-braid of length $\ell$ where the closure is a knot satisfies 
    \begin{equation*}
        \lim_{\ell\rightarrow\infty} \frac{E(v_2)}{\ell} = \frac{1}{8}.
    \end{equation*}
\end{theorem}

Thus for the crossing number $c$ of such a knot we have     $\displaystyle \lim_{c\rightarrow\infty} \frac{E(v_2)}{c} \geq 
\frac{1}{8}$.  
These knots have bridge number at most 3.  Comparing our Main Theorem \ref{thm:main} to Theorem \ref{thm:3braid} by \cite{Ville:Casson}, perhaps we may speculate that the random braids can be simplified some.

\section{Background}
\label{sec:background}

We begin with one of the results used in the proof of the main Theorem \ref{thm:main} above.  Ernst and Sumners computed the number of 2-bridge knots with chiral pairs \textit{not} counted separately.
\begin{theorem}[Ernst-Sumners \cite{ErnSum}, Theorem 5]
\label{thm:ernstsumners}
The number $|\mathcal{K}(c)|$ of 2-bridge knots with $c$ crossings where chiral pairs are \emph{not} counted separately is given by
\[
|\mathcal{K}(c)| = 
\begin{cases}
\frac{1}{3}(2^{c-3}+2^{\frac{c-4}{2}}) & \text{ for }4 \leq c\equiv 0 \text{ mod }4,\\
\frac{1}{3}(2^{c-3}+2^{\frac{c-3}{2}}) & \text{ for }5\leq c\equiv 1 \text{ mod }4, \\
\frac{1}{3}(2^{c-3}+2^{\frac{c-4}{2}}-1) & \text{ for }6 \leq c\equiv 2 \text{ mod }4, \text{ and}\\
\frac{1}{3}(2^{c-3}+2^{\frac{c-3}{2}}+1) & \text{ for }3\leq c\equiv 3 \text{ mod }4.
\end{cases}
\]
\end{theorem}

We will consider twice this set $\mathcal{K}(c)$ and ``partition'' it into two sets $T(c)$ and $T_p(c)$.  These were introduced by the first author with Lowrance in \cite{CohLow} based on work with Krishnan and Even-Zohar in \cite{CoKr, CoEZKr, Coh:lower}.  In turn, these articles were based on work by Koseleff and Pecker with others \cite{KosPec1, KosPecRou:10crossings, KosPec4, KosPec3,  KosPec:billiard, KosPec:Alex, KosPec, KosPecTran:diagrams, BruKosPec}.

\begin{definition}
    \label{def:tc}
The \textbf{partially doubly-counted set $T(c)$ of 2-bridge knots with crossing number $c$} is the set of words in $\{+,-\}$ of the form
$$(+)^{\varepsilon_1}(-)^{\varepsilon_2}(+)^{\varepsilon_3}(-)^{\varepsilon_4}\ldots(-)^{\varepsilon_{c-1}}(+)^{\varepsilon_{c}}$$
for odd $c$ or of the form 
$$(+)^{\varepsilon_1}(-)^{\varepsilon_2}(+)^{\varepsilon_3}(-)^{\varepsilon_4}\ldots(+)^{\varepsilon_{c-1}}(-)^{\varepsilon_{c}}$$
for even $c$, where in both cases $\varepsilon_1=\varepsilon_c=1$, for $i\in\{2,\ldots,c-1\}$ the remaining $\varepsilon_i\in\{1,2\}$, and the (reduced) length of the word $\ell=\sum_{i=1}^c \varepsilon_i\equiv 1$ mod 3.

Such a knot is \textbf{palindromic} if $\varepsilon_i=\varepsilon_{c+1-i}$ for all $1 \leq  i \leq c$.

Let $T_p(c)$ be the \textbf{set of palindromic 2-bridge knots with crossing number $c$}.  This is a subset of $T(c)$.

Each of the words $w$ in $T(c)$ or $T_p(c)$ corresponds to an alternating 2-bridge knot diagram $K_w$ with $(+)$ and $(-)(-)$ corresponding to crossings \crosspos \; or $\sigma_1$ appearing in the lower row and $(-)$ and $(+)(+)$ corresponding to crossings \crossneg \; or $\sigma_2^{-1}$ appearing in the upper row.  These braids do not have the usual braid closure but instead have a plat closure where the 3-braid is taken with a 4th unaffected strand assumed to be above.

Observe that the final crossing for odd $c$ is $\sigma_1$ and that the final crossing for even $c$ is $\sigma_2^{-1}$.
\end{definition}

Palindromic 2-bridge knots appear just once in $T(c)$; every other word in $T(c)$ appears both forward and backward, and hence the name of this set.

\begin{figure}[h!]
\begin{tikzpicture}
\draw[thick, -stealth] (-1,0) -- (0,0);
\draw[thick, -stealth] (6,2) -- (7,2);
\draw[thick, -stealth] (0,2) arc (90:270:.5);
\draw[thick, stealth-] (6,0) arc (-90:90:.5);
\draw[thick, stealth-] (0,2) -- (1,2);
\draw[thick, -stealth] (1,0) -- (2,0);
\draw[thick, -stealth] (2,2) -- (5,2);
\draw[thick, stealth-] (5,0) -- (6,0);
\foreach \x/ \y in {0/0}
    {
    \draw[thick, stealth-] (\x+1,\y) -- (\x,\y+1);
    \draw[color=white, line width=10] (\x,\y) -- (\x+1,\y+1);
    \draw[thick, -stealth] (\x,\y) -- (\x+1,\y+1);
    }
\foreach \x/ \y in {5/1}
    {
    \draw[thick, -stealth] (\x,\y) -- (\x+1,\y+1);
    \draw[color=white, line width=10] (\x+1,\y) -- (\x,\y+1);
    \draw[thick, stealth-] (\x+1,\y) -- (\x,\y+1);
    }
\foreach \x/ \y in {1/1, 2/0, 4/0}
    {
    \draw[thick, -stealth] (\x,\y) -- (\x+1,\y+1);
    \draw[color=white, line width=10] (\x+1,\y) -- (\x,\y+1);
    \draw[thick, -stealth] (\x+1,\y) -- (\x,\y+1);
    }
\foreach \x/ \y in {2/0, 4/0}
    {
    \draw[thick, -stealth] (\x+1,\y) -- (\x,\y+1);
    \draw[color=white, line width=10] (\x,\y) -- (\x+1,\y+1);
    \draw[thick, -stealth] (\x,\y) -- (\x+1,\y+1);
    }
\foreach \x/ \y in {3/0}
    {
    \draw[thick, stealth-] (\x+1,\y) -- (\x,\y+1);
    \draw[color=white, line width=10] (\x,\y) -- (\x+1,\y+1);
    \draw[thick, stealth-] (\x,\y) -- (\x+1,\y+1);
    }
\filldraw[black] (-1,0) circle (2pt);
\end{tikzpicture}
    \caption{\label{fig:61} The stevedore knot 
    as in Example \ref{ex:61}}
\end{figure}
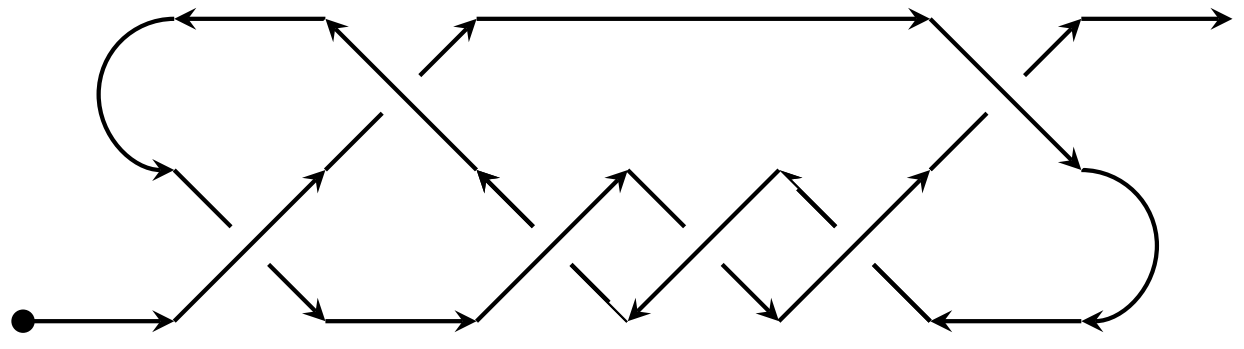

\begin{example}
    \label{ex:61}
    The stevedore knot $6_1$ can be written $\sigma_1\sigma_2^{-1}\sigma_1^3\sigma_2^{-1}$ corresponding to the word $+-+-{}-+-$ of length 7 in the set $T(6)$ as in Figure \ref{fig:61}.
\end{example}

We have that $2|\mathcal{K}(c)|=|T(c)|+|T_p(c)|$ and abuse notation to say that the latter two sets partition twice the first set.

The first author with Lowrance \cite[Proposition 3.1]{CohLow} proved that  $|T(c)|=t(c)=\frac{1}{3}(2^{c-2}(-1)^c)$ is the $(c-2)$nd Jacobsthal number $J(c-2)$ that satisfies the Fibonacci-like recursion 
\begin{equation}
    \label{eq:tc}
    t(c)=t(c-1)+2t(c-2).
\end{equation}
This is proved by way of a bijection between the set $T(c)$ and the set $T(c-1)\sqcup T(c-2) \sqcup T(c-2)$, yielding a recursion that will appear in Subsection \ref{subsec:recursion} and Figure \ref{fig:cases}.

The set $T(c)$ is partitioned into four subsets $T_i(c)$ of size $t_i(c)$ corresponding to each of the four cases below, 
and it is shown \cite[Lemma 3.4]{CohLow} that 
\begin{equation}
\label{eq:tic}
t_1(c)=2t(c-3) \text{ and }t_2(c)=t_3(c)=t_4(c)=t(c-2).
\end{equation}

This bijection is achieved by taking a word in $T(c)$, removing the last three crossings, which we call \textbf{tails}, and replacing them with one or two crossings, the \textbf{new tails}, based on four cases (described here for odd $c$):  
\begin{itemize}
    \item[(I)]  $w\sigma_1\sigma_2^{-1}\sigma_1$ is replaced by  $w\sigma_2^{-1}\sigma_2^{-1}$ in $T(c-1)$,
    \item[(II)] $w\sigma_2^{-1}\sigma_1\sigma_1$ is replaced by $w\sigma_1\sigma_2^{-1}$ in $T(c-1)$,
    \item[(III)] $w\sigma_1\sigma_1\sigma_1$ is replaced by $w\sigma_1$ in $T(c-2)$, and
\item[(IV)] $w\sigma_2^{-1}\sigma_2^{-1}\sigma_1$ is replaced by $w\sigma_1$ in $T(c-2)$.
\end{itemize}
It is important to note that these replacements preserve the subword $w$ as well as the components of the tails, and so the orientations of all of the previous crossings remain unchanged.  Also note that cases I and II partition $T(c-1)$.

We now make the following definitions with the goal of re-interpreting the Casson invariant in a more accessible way as in Polyak and Viro \cite{PolVir:Vas}.

\begin{definition}
    \label{def:pairs}
    We consider each knot as a long knot with the marked point at infinity drawn at the lower left.  Starting at this marked point, we travel along the knot following the orientation.  We denote a crossing by the pair $(p,\tilde{p})$ where $p$ is the first time we arrive at the crossing and $\tilde{p}$ is the second time.  
    
The \textbf{sign of a crossing} of an oriented knot is given by $sign(\crossposor)=+1$ and $sign(\crossnegor)=-1$.

    Let $i(p)=i((p,\tilde{p}))=i(\tilde{p})$ be the \textbf{index} of the crossing $p$, that is, the order of its crossing from left to right between $1\leq i(p) \leq c$.
\end{definition}

\begin{definition}
\label{def:ascdesc}
A crossing $(p,\tilde{p})$ is \textbf{ascending} if $p$ belongs to the understrand of the crossing or \textbf{descending} if $p$ belongs to the overstrand of the crossing.

A pair of crossings $((p,\tilde{p}),(q,\tilde{q}))$ is \textbf{valid} if 
\begin{enumerate}
	\item the crossing $p$ is ascending,
	\item the crossing $q$ is descending, and 
	\item the crossings are alternating:  that is, they  appear in the order $p,q,\tilde{p},\tilde{q}$ along the knot starting at the marked point. 
\end{enumerate}
\end{definition}

\begin{figure}[h!]
\begin{tikzpicture}
\draw[thick, -stealth] (-1,0) -- (0,0);
\draw[thick, -stealth] (6,2) -- (7,2);
\draw[thick, -stealth] (0,2) arc (90:270:.5);
\draw[thick, stealth-] (6,0) arc (-90:90:.5);
\draw[thick, stealth-] (0,2) -- (1,2);
\draw[thick, -stealth] (1,0) -- (2,0);
\draw[thick, -stealth] (2,2) -- (5,2);
\draw[thick, stealth-] (5,0) -- (6,0);
\foreach \x/ \y in {0/0}
    {
    \draw[thick, stealth-] (\x+1,\y) -- (\x,\y+1);
    \draw[color=white, line width=10] (\x,\y) -- (\x+1,\y+1);
    \draw[thick, -stealth] (\x,\y) -- (\x+1,\y+1);
    }
\foreach \x/ \y in {5/1}
    {
    \draw[thick, -stealth] (\x,\y) -- (\x+1,\y+1);
    \draw[color=white, line width=10] (\x+1,\y) -- (\x,\y+1);
    \draw[thick, stealth-] (\x+1,\y) -- (\x,\y+1);
    }
\foreach \x/ \y in {1/1, 2/0, 4/0}
    {
    \draw[thick, -stealth] (\x,\y) -- (\x+1,\y+1);
    \draw[color=white, line width=10] (\x+1,\y) -- (\x,\y+1);
    \draw[thick, -stealth] (\x+1,\y) -- (\x,\y+1);
    }
\foreach \x/ \y in {2/0, 4/0}
    {
    \draw[thick, -stealth] (\x+1,\y) -- (\x,\y+1);
    \draw[color=white, line width=10] (\x,\y) -- (\x+1,\y+1);
    \draw[thick, -stealth] (\x,\y) -- (\x+1,\y+1);
    }
\foreach \x/ \y in {3/0}
    {
    \draw[thick, stealth-] (\x+1,\y) -- (\x,\y+1);
    \draw[color=white, line width=10] (\x,\y) -- (\x+1,\y+1);
    \draw[thick, stealth-] (\x,\y) -- (\x+1,\y+1);
    }
\filldraw[black] (-1,0) circle (2pt);

\node[below] at (.5,0) {$q_1^-$};
\node[above] at (1.5,2) {$p_2^-$};
\node[below] at (2.5,0) {$p_3^+$};
\node[below] at (3.5,0) {$q_4^+$};
\node[below] at (4.5,0) {$p_5^+$};
\node[above] at (5.5,2) {$q_6^+$};

\end{tikzpicture}
    \caption{\label{fig:61letters} The stevedore knot $6_1$ from Figure \ref{fig:61} with its six crossings labeled ascending ($p_i$) or descending ($q_i$) and positive or negative as in Definitions \ref{def:pairs} and \ref{def:ascdesc} 
    as in Example \ref{ex:61crossings}} 
\end{figure}
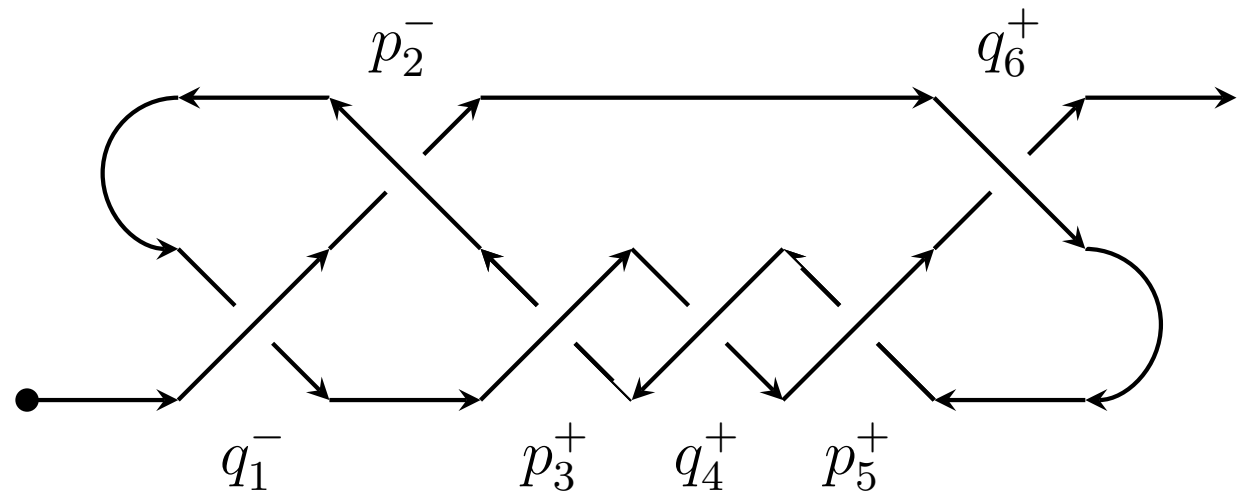

\begin{notation}
\label{not:pq}
        In Figure \ref{fig:61letters} and others below, we refer to crossings as $p_i^\pm$ and $q_i^\pm$ where $p$ or, respectively, $q$ denote ascending or, respectively, descending, the $i$ denotes the index of the crossing, and $\pm$ denotes the sign of the crossing.
\end{notation}

We note that the signs of the crossings for diagrams coming from words in $T(c)$ can be easily determined through work of the first author \cite[Property 7.1]{Co:3-bridge} and \cite{CohLow}.

\begin{example}
    \label{ex:61crossings}
    The first two crossings of the stevedore knot $6_1$ are negatively signed and the rest are positively signed.  Three crossings are ascending and three are descending as seen in Figure \ref{fig:61}.  

The pair ($p_3^+,q_4^+$) is not a valid pair because $q_4$ is reached before $p_3$ when traveling along the knot from the marked point.  Similarly, the pairs ($p_3^+,q_6^+$) and ($p_5^+,q_6^+$) are not valid.  The only two valid pairs are ($p_2^-,q_4^+$) and ($p_2^-,q_6^+$)  because they alternate correctly.  
\end{example}

\begin{remark}
    In the next section, we give a criterion to describe how to more easily find these valid pairs.  See Lemma \ref{lem:validcrossings} and Example \ref{ex:pairs}.
\end{remark}

The following result comes from Polyak and Viro \cite{PolVir:Vas} but is more explicitly stated in this form by Even-Zohar, Hass, Linial, and Nowik \cite[Lemma 8]{EZHLN}.

\begin{proposition}[Polyak-Viro \cite{PolVir:Vas}]
    \label{prop:cassonpairs}
The Casson invariant is
      \begin{equation}
          \label{eq:vdef}
            v_2(K)=\sum sign(p)sign(q)
      \end{equation}
summing over all valid pairs of crossings.
\end{proposition}

The computations of the Casson invariants of several knots using this result appear in 
Subsection \ref{subsec:examples} below after some useful preliminary results.

\section{Contribution from the partially doubly-counted set}
\label{sec:main}

Recall from the introduction 
that $v(c)=\sum_{K\in T(c)} v_2(K)$ is the sum of the Casson invariant of each knot in the partially doubly-counted set $T(c)$.  The goal of this section is to prove the following theorem, which is used in the proof of Main Theorem \ref{thm:main} above.

\begin{theorem}
    \label{thm:v}
The total contribution to the Casson invariant over all knots from the set $T(c)$ has the following closed formula:
\begin{equation}
\label{eq:mainv}
v(c) \;  = \; \frac{1}{72}c2^{c}  +\frac{1}{432}2^{c}   -\frac{3}{16}(-2)^{c} + \frac{2}{9}(-1)^cc   +\frac{5}{27}(-1)^{c}.
\end{equation}
\end{theorem}

Recall Definitions \ref{def:pairs} and \ref{def:ascdesc}, and recall from Proposition \ref{prop:cassonpairs} the formula for the Casson invariant $v_2(K)$ summing over \emph{valid} pairs $p$ and $q$ of crossings with $p$ ascending and with $q$ descending.  

\begin{definition}
    \label{def:strands}
We define the \textbf{first} \textcircled{1}, \textbf{second} \textcircled{2}, and \textbf{third} \textcircled{3} \textbf{strands} to be the portions of the knot diagram traversed from one side of the knot to the other side in order starting from the marked point.
\end{definition}

\begin{figure}[h!]
\begin{tikzpicture}
\draw[ultra thick, -stealth] (-1,0) -- (0,0) -- (1,1) -- (1.2,1.2);
\draw[ultra thick, -stealth] (1.8,1.8) -- (2,2) -- (5,2) --  (6,1);

\draw[dashed, -stealth] (6,0) -- (5,0) -- (4.8,.2);
\draw[dashed, -stealth] (4.2,.8) -- (4,1) -- (3,0) -- (2.8, .2);
\draw[dashed, -stealth] (2.2,.8) -- (2,1) -- (1,2) -- (0,2);

\draw[thick, -stealth] (0,1) -- (.2,.8);
\draw[thick, -stealth] (.8,.2) -- (1,0) -- (2,0) -- (3,1) -- (3.2,.8);
\draw[thick, -stealth] (3.8,.2) -- (4,0) -- (5.2,1.2);
\draw[thick, -stealth] (5.8,1.8) -- (6,2) -- (7,2);

\draw[-stealth] (0,2) arc (90:270:.5);
\draw[stealth-] (6,0) arc (-90:90:.5);

\node[below] at (-.5,0) {\textcircled{1}};
\node[below] at (5.5,0) {\textcircled{2}};
\node[above] at (0,1) {\textcircled{3}};

\filldraw[black] (-1,0) circle (3pt);
\end{tikzpicture}
    \caption{\label{fig:61firstsecondthird} The stevedore knot $6_1$ from Figure \ref{fig:61}  with its three strands labeled following the orientation from the marked point as in Definition \ref{def:pairs}:  here the first strand is the thickest and the second strand is dashed}
\end{figure}
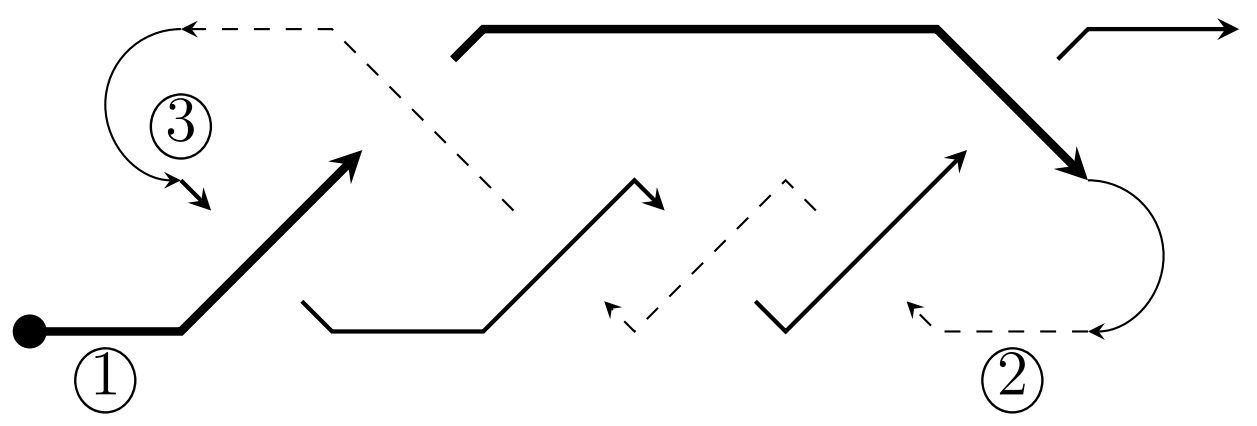

To illustrate the numbered strands see Figure \ref{fig:61firstsecondthird}.

We begin with a result which characterizes valid pairs of crossings.

\begin{lemma}
\label{lem:validcrossings}
    The following is true for all \emph{valid} pairs of crossings $(p,\tilde{p})$ and $(q,\tilde{q})$.  
    
    The index $i(p)$ of the ascending crossing must be less than the index $i(q)$ of the descending crossing.
    
    The (first appearance of the) ascending crossing $p$ must lie on the first strand \emph{\textcircled{1}}.  The (second appearance of the) descending crossing $\tilde{q}$ must lie on the third strand \emph{\textcircled{3}}.  The only four possible cases are when
    \begin{enumerate}
        \item  $\tilde{p}$ lies on the second strand \emph{\textcircled{2}} and $q$ on the first \emph{\textcircled{1}},
        \item  $\tilde{p}$ lies on the third strand \emph{\textcircled{3}} and $q$ on the first \emph{\textcircled{1}},
        \item  $\tilde{p}$ lies on the second strand \emph{\textcircled{2}} and $q$ on the second \emph{\textcircled{2}}, or 
        \item  $\tilde{p}$ lies on the third strand \emph{\textcircled{3}} and $q$ on the second \emph{\textcircled{2}}.      
    \end{enumerate}
\end{lemma}

\begin{proof} Note that $p$ and $\tilde{p}$ (resp. $q$ and $\tilde{q}$) cannot belong to the same strand.

    Suppose first that $p$ and $q$ are in \textcircled{1}.  Then $i(p)< i(q)$. 
    If $\tilde{p}$ and $\tilde{q}$ are in \textcircled{2}, then $\tilde{q}$ comes before $\tilde{p}$, which is not allowed for a valid pair: \textcircled{2} goes from right to left, whereas \textcircled{1} goes from left to right.  Then a valid pair must be in one of the following cases:  either (1) $\tilde{p}\in $\textcircled{2} and $\tilde{q}\in $\textcircled{3} or (2) both $\tilde{p}, \tilde{q}\in$\textcircled{3}.
    
Suppose next that $p$ is in \textcircled{1} but that $q$  in \textcircled{2}. If (3) $\tilde{p}\in $\textcircled{2} and $\tilde{q}\in$ \textcircled{3}, then since $q$ comes before $\tilde{p}$, we have $i(q)>i(\tilde{p})$, so $i(p)<i(q)$.  If (4) both $\tilde{p}, \tilde{q}\in$\textcircled{3}, then $i(\tilde{p})<i(\tilde{q})$ or $i(p)<i(q)$.

Finally suppose that both $p, q$ are in \textcircled{2}.  Then  $i(q)<i(p)$.  But since both $\tilde{p}, \tilde{q}$ are in \textcircled{3}, they satisfy  $i(\tilde{q})< i(\tilde{p})$, a contradiction.
\end{proof}

\begin{example}
    \label{ex:pairs} In Figures \ref{fig:61letters} and \ref{fig:61firstsecondthird} for the stevedore knot $6_1$, we see that the only ascending crossing that is on strand \textcircled{1} is $p_2$, so neither $p_3$ nor $p_5$ can be a part of a valid pair.  Furthermore, we need $i(p_2)<i(q)$, and so $q_1$ cannot be a part of a valid pair.  This makes it easier to see that the only valid pairs are ($p_2^-,q_4^+$) and ($p_2^-,q_6^+$), as in Example \ref{ex:61crossings}.
\end{example}

In order to perform the replacements of Section \ref{sec:background}, we now show where the first and third strands lie in the tail.

\begin{lemma}\label{lem:penultimateorder} Between the penultimate and last crossings of a knot $K_c\in T(c)$, the strands are in the following order going from top to bottom if $c$ is even (respectively odd):
$$\begin{pmatrix}1\\ 3\\ 2\end{pmatrix}\ \ \ \left(\mbox{resp.}\ \begin{pmatrix}2\\ 3\\  1\end{pmatrix}\right).$$
\end{lemma}

\begin{proof}
    This follows from the figures in Table 4 of Cohen-Lowrance \cite{CohLow}, which are recreated below in Figure \ref{fig:cases}, and the arguments that the second strand is the only one oriented to the left and the third strand is the long strand at the right end.
\end{proof}

To more easily perform calculations giving us a closed formula for $v(c)$, we  define a new diagrammatic signed counting of the ascending crossings $(p,\tilde{p})$ on the first strand. 

\begin{definition}
    \label{def:f}
Let $w$ be a word in $T(c)$, and let $D_w$ be its associated alternating knot diagram.  

    We define the \textbf{first-strand signed ascending sum} of a word (or knot diagram) to be $f(w)=f(D_w)=\displaystyle\sum_{p\text{ ascending on \textcircled{1}}} sign(p)$ and for a given crossing number to be\\ $f(c)=\displaystyle\sum_{w\in T(c)} f(w)$.

    Note that the first and final crossings are always descending, so they will never be counted by the first-strand signed ascending summation.
    \end{definition}

\subsection{Examples for small crossing numbers} 
\label{subsec:examples} 
We will use these knots as base cases in the inductive proofs below.

\begin{figure}[ht]
\begin{tikzpicture}
\node at (-1.5,2){\textcircled{2}};
\node at (-1.5,){\textcircled{3}};
\node at (-1.5,0){\textcircled{1}};
\draw[thick, -stealth] (-1,0) -- (0,0);
\draw[thick, -stealth] (5,0) -- (6,0);
\draw[thick, -stealth] (0,2) arc (90:270:.5);
\draw[thick, -stealth] (5,1) arc (-90:90:.5);
\draw[thick, stealth-] (0,2) -- (5,2);
\foreach \x/ \y in {0/0, 2/0, 4/0}{
\draw[thick, stealth-] (\x+1,\y) -- (\x,\y+1);
\draw[color=white, line width=10] (\x,\y) -- (\x+1,\y+1);
\draw[thick, -stealth] (\x,\y) -- (\x+1,\y+1);
}

\filldraw[black] (-1,0) circle (2pt);
    \node at (0.5,-0.4)  {$q_1^-$};
\node at (2.5,-0.4)  {$p_2^-$};
\node at (4.5,-0.4)  {$q_3^-$};

\draw[thick, -stealth] (1,0) -- (2,0);
\draw[thick, -stealth] (1,1) -- (2,1);
\draw[thick, -stealth] (3,0) -- (4,0);
\draw[thick, -stealth] (3,1) -- (4,1);

\end{tikzpicture}
    \caption{\label{fig:3} The trefoil knot $3_1$ with strands and crossings labeled }
\end{figure}
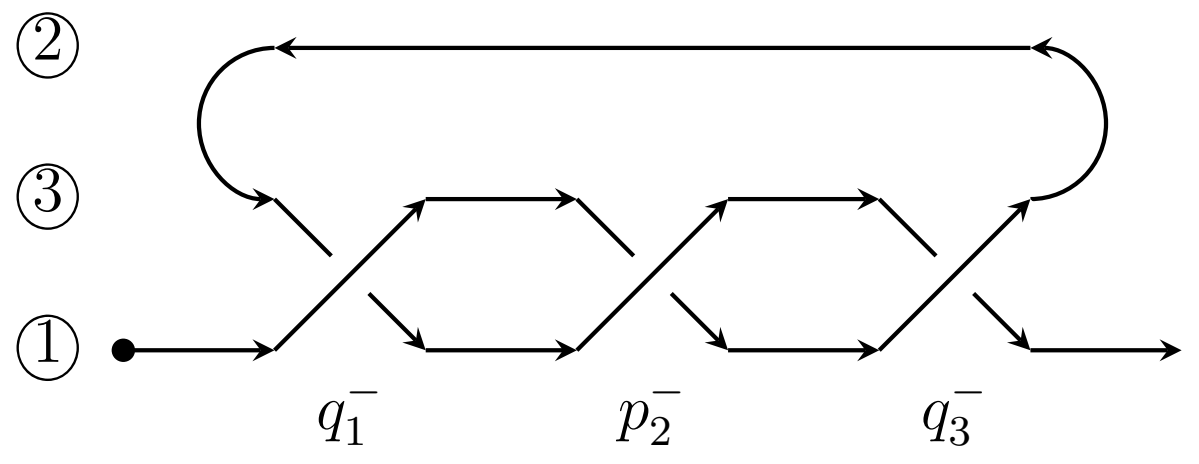

There is only one knot in $T(3)$:  the trefoil $3_1$ written $+--+$ or $\sigma_1^3$ with all crossings negatively signed, first descending, then ascending, and then descending.  The only valid pair is ($p_2^-,q_3^-$) so $v(3)=1$.  The only ascending crossing on the first strand \textcircled{1} is $p_2^-$ so we have  $f(3)=-1$.  See Figure \ref{fig:3}.

\begin{figure}[ht]
      \begin{tikzpicture}
            \node at (-1.5,2){\textcircled{2}};
      \node at (-1.5,){\textcircled{3}};
      \node at (-1.5,0){\textcircled{1}};
      \draw[thick, -stealth] (-1,0) -- (0,0);
      \draw[thick, -stealth] (4,0) -- (3,0);
      \draw[thick, -stealth] (0,2) arc (90:270:.5);
      \draw[thick, stealth-] (4,0) arc (-90:90:.5);
      \draw[thick, stealth-] (0,2) -- (1,2);
      \draw[thick, -stealth] (1,0) -- (2,0);
      \draw[thick, -stealth] (2,2) -- (3,2);
      \draw[thick, -stealth] (4,2) -- (5,2);
\foreach \x/ \y in {0/0}{
      \draw[thick, stealth-] (\x+1,\y) -- (\x,\y+1);
      \draw[color=white, line width=10] (\x,\y) -- (\x+1,\y+1);
      \draw[thick, -stealth] (\x,\y) -- (\x+1,\y+1);  
      }
\foreach \x/ \y in {1/1}{
    \draw[thick, -stealth] (\x,\y) -- (\x+1,\y+1);
    \draw[color=white, line width=10] (\x+1,\y) -- (\x,\y+1);
    \draw[thick, -stealth] (\x+1,\y) -- (\x,\y+1);
      }
\foreach \x/ \y in {2/0}{
      \draw[thick, -stealth] (\x+1,\y) -- (\x,\y+1);
      \draw[color=white, line width=10] (\x,\y) -- (\x+1,\y+1);
      \draw[thick, -stealth] (\x,\y) -- (\x+1,\y+1);  
      }
\foreach \x/ \y in {3/1}{
    \draw[thick, -stealth] (\x,\y) -- (\x+1,\y+1);
    \draw[color=white, line width=10] (\x+1,\y) -- (\x,\y+1);
    \draw[thick, stealth-] (\x+1,\y) -- (\x,\y+1);
      }

      \filldraw[black] (-1,0) circle (2pt);
    \node at (0.5,-0.4)  {$q_1^-$};
      \node at (1.5,2.4)  {$p_2^-$};
      \node at (2.5,-0.4)  {$p_3^+$};
      \node at (3.5,2.4)  {$q_4^+$};
      
\end{tikzpicture}
    \caption{\label{fig:4} The figure eight knot $4_1$ with strands and crossings labeled }
\end{figure}
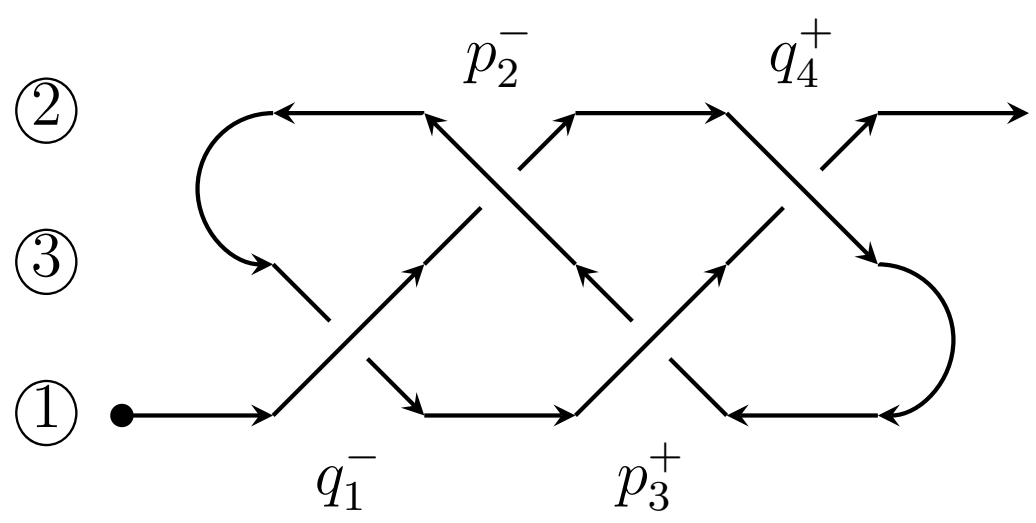

There is only one knot in $T(4)$:  the figure eight $4_1$ written $+-+-$ or $\sigma_1\sigma_2^{-1}\sigma_1\sigma_2^{-1}$ with  the first two crossings negatively signed, first descending then ascending, and the last two crossings positively signed, first ascending then descending.  The only valid pair is ($p_2^-,q_4^+$) so $v(4)=-1$.  The only ascending crossing on the first strand \textcircled{1} is $p_2^-$ so  $f(4)=-1$.  
  See Figure \ref{fig:4}.

\begin{figure}[ht]
            \begin{tikzpicture}
            \node at (-1.5,2){\textcircled{2}};
            \node at (-1.5,){\textcircled{3}};
            \node at (-1.5,0){\textcircled{1}};
            \draw[thick, -stealth] (-1,0) -- (0,0);
            \draw[thick, -stealth] (2,0) -- (4,0);
            \draw[thick, stealth-] (5,1) -- (4,0);
            \draw[thick, -stealth] (5,1) -- (6,1);
            \draw[thick, -stealth] (0,2) arc (90:270:.5);
            \draw[thick, stealth-] (7,1) arc (-90:90:.5);
            \draw[thick, stealth-] (0,2) -- (1,2);
            \draw[thick, -] (1,0) -- (3,0);
            \draw[thick, -stealth] (2,2) -- (3,2);
            \foreach \x/ \y in {0/0}
            {
                  \draw[thick, stealth-] (\x+1,\y) -- (\x,\y+1);
                  \draw[color=white, line width=10] (\x,\y) -- (\x+1,\y+1);
                  \draw[thick, -stealth] (\x,\y) -- (\x+1,\y+1);  
            }
            \draw[thick, -stealth] (6,1) -- (6.4,0.6);
            \draw[thick, -stealth] (6.6,0.4) -- (7,0);
            \draw[thick, stealth-] (1,2) -- (2,1);
            \draw[thick, -stealth] (3,1) -- (2,1);
            \draw[thick, -stealth] (3.4,1.4) -- (3,1);
            \draw[thick, stealth-] (3.6,1.6) -- (4,2);
            \draw[thick, -stealth] (1.6,1.6) -- (2,2);
            \draw[thick] (1,1) -- (1.3,1.3);
            \draw[thick, -stealth] (3,2) -- (4,1) -- (4.4,0.6);
            \draw[thick, -stealth] (4.6,0.4) -- (5,0);
            \draw[thick, stealth-] (6,0) -- (5,0);
            \draw[thick, -stealth] (6,0) -- (7,1);
            \draw[thick, stealth-] (4,2) -- (7,2);
            \node at (0.5,-0.3){$q_1^-$};
            \node at (1.5,2.4) {$p_2^-$};
            \node at (3.5,2.4) {$q_3^-$};
            \node at (4.5,-0.3){$p_4^-$};
            \node at (6.5,-0.3){$q_5^-$};
            
            \filldraw[black] (-1,0) circle (2pt);
            \end{tikzpicture}

    \begin{tikzpicture}
            \node at (-1.5,2){\textcircled{2}};
      \node at (-1.5,){\textcircled{3}};
      \node at (-1.5,0){\textcircled{1}};
      \draw[thick, -stealth] (-1,0) -- (0,0);

      \draw[thick, -stealth] (0,2) arc (90:270:.5);
      \draw[thick, stealth-] (7,1) arc (-90:90:.5);
      \draw[thick, -stealth] (1,1) -- (2,1);
      
      \draw[thick, stealth-] (0,2) -- (3,2);
      \foreach \x/ \y in {0/0}
      {
            \draw[thick, stealth-] (\x+1,\y) -- (\x,\y+1);
            \draw[color=white, line width=10] (\x,\y) -- (\x+1,\y+1);
            \draw[thick, -stealth] (\x,\y) -- (\x+1,\y+1);  
      }
      \draw[thick, -stealth] (6,1) -- (6.4,0.6);
      \draw[thick, -stealth] (6.6,0.4) -- (7,0);

\draw[thick, stealth-] (6,2) -- (7,2);
      
      \draw[thick, -stealth] (3.6,1.6) -- (4,2);
      
      \draw[thick, stealth-] (3,2) -- (4,1);
      \draw[thick, -stealth] (5,1) -- (4,1);

      \draw[thick, -stealth] (3,0) -- (6,0);
      \draw[thick, -stealth] (6,0) -- (7,1);
      \draw[thick,- stealth] (4,2) -- (5,2);
      \draw[thick,- stealth] (5.4,1.4) -- (5,1);
      \draw[thick,stealth-] (5.6,1.6) -- (6,2);
      \draw[thick,stealth-] (6,1) -- (5,2);
      \draw[thick, stealth-] (3,0) -- (2.6,0.4);
      \node at (0.5,-0.3){$q_1^-$};
      \node at (2.5,-0.3) {$p_2^-$};
      \node at (3.5,2.4)  {$q_3^-$};
      \node at (5.3,2.4){$p_4^-$};
      \node at (6.5,-0.3){$q_5^-$};
            \draw[thick, -stealth] (2,0) -- (3.4,1.4);
      \draw[thick] (2,1) -- (2.4,0.6);
      
      \draw[thick] (1,0) -- (2,0);

      \filldraw[black] (-1,0) circle (2pt);
      \end{tikzpicture}
    \caption{\label{fig:52} Two copies of the 3-twist knot $5_2$ with strands and crossings labeled}
\end{figure}
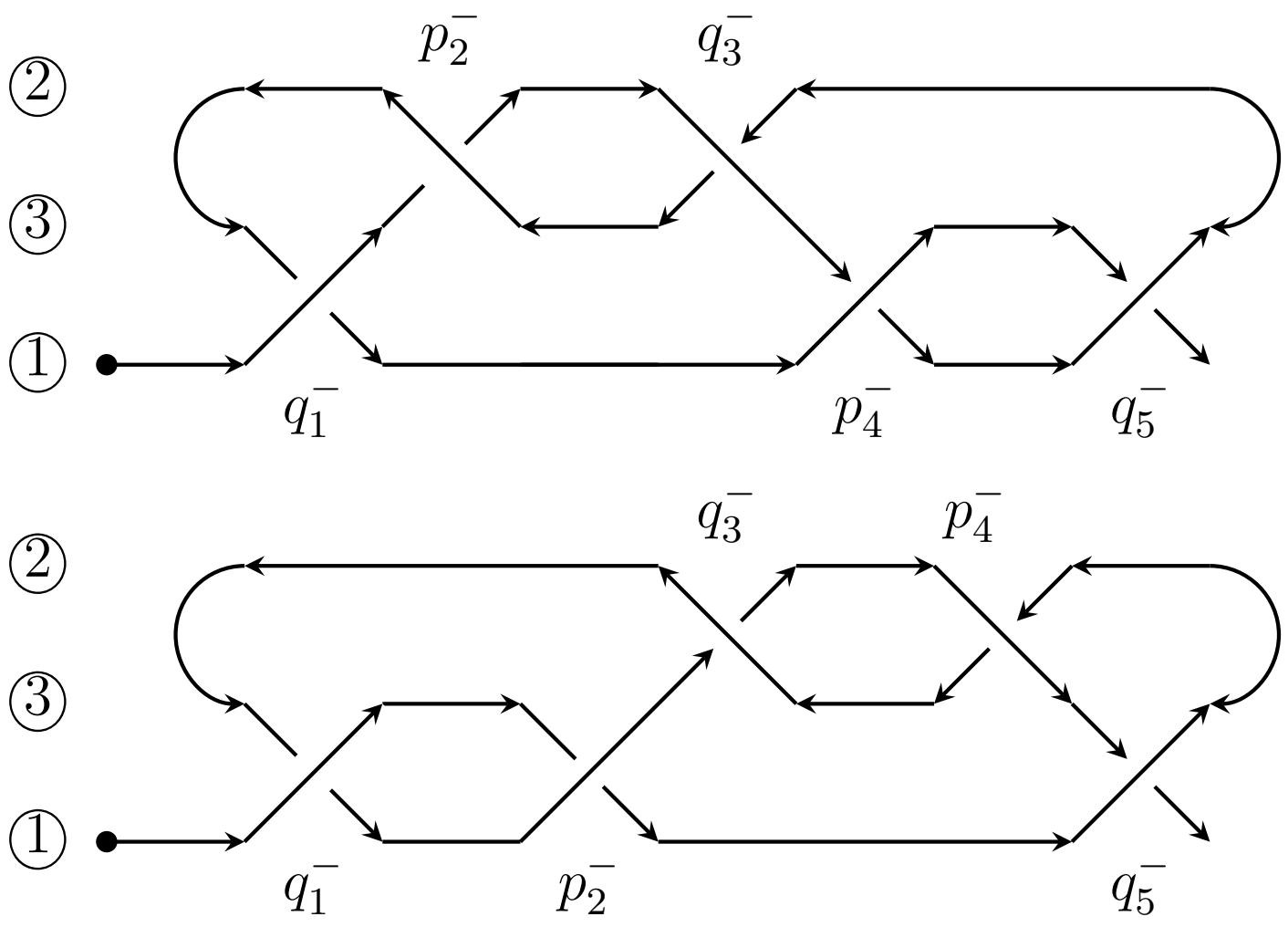

There are three knots in $T(5)$:  the pentafoil $5_1$ or torus knot $T(2,5)$ written $+-{}-+-{}-+$ or $\sigma_1^5$ and the 3-twist knot $5_2$ counted twice written $+-++-{}-+$ or $\sigma_1\sigma_2^{-2}\sigma_1^2$ along with its reverse.  The first has three valid pairs ($p_2^-,q_3^-$), ($p_2^-,q_5^-$), and ($p_4^-,q_5^-$), the second has two valid pairs ($p_2^-,q_5^-$) and ($p_4^-,q_5^-$), and the third has two valid pairs 
($p_2^-,q_3^-$) and ($p_2^-,q_5^-$),
 so $v(5)=[3+2+2]=7$.  We also have that $f(5)=[(-2)+(-2)+(-1)]=-5$. 
See Figures \ref{fig:51} and \ref{fig:52}.  Compare the two diagrams in Figure \ref{fig:52} to see that $f$ is only a diagrammatic count and is not an invariant of the knot type.

      \begin{figure}[ht]
      \begin{tikzpicture}     
            \node at (-1.5,2){\textcircled{2}};
      \node at (-1.5,){\textcircled{3}};
      \node at (-1.5,0){\textcircled{1}};
      \draw[thick, -stealth] (0,2) arc (90:270:.5);
\draw[thick, -stealth] (9,1) arc (-90:90:.5);
      \draw[thick, -stealth] (0,0) -- (1,1);
      \draw[thick] (0,1) -- (0.4,0.6);
      \draw[thick, stealth-] (1,0) -- (0.6,0.4);
      \draw[thick] (-1,0) -- (0,0);

      \draw[thick, -stealth] (2,0) -- (3,1);
      \draw[thick] (2,1) -- (2.4,0.6);
      \draw[thick, stealth-] (3,0) -- (2.6,0.4);
      \draw[thick] (1,0) -- (2,0);

      \draw[thick, -stealth] (4,0) -- (5,1);
      \draw[thick] (4,1) -- (4.4,0.6);
      \draw[thick, stealth-] (5,0) -- (4.6,0.4);
      \draw[thick] (3,0) -- (4,0);

\draw[thick, -stealth] (6,0) -- (7,1);
\draw[thick] (6,1) -- (6.4,0.6);
\draw[thick, stealth-] (7,0) -- (6.6,0.4);
\draw[thick] (5,0) -- (6,0);

\draw[thick, -stealth] (8,0) -- (9,1);
\draw[thick] (8,1) -- (8.4,0.6);
\draw[thick, stealth-] (9,0) -- (8.6,0.4);
\draw[thick] (7,0) -- (8,0);

\draw[thick] (1,1) -- (2,1);
\draw[thick] (3,1) -- (4,1);

\draw[thick] (5,1) -- (6,1);
\draw[thick] (7,1) -- (8,1);
\draw[thick, -stealth] (9,0) -- (10,0);
\draw[thick, stealth-] (0,2) -- (9,2);
\node at (0.5,-0.3){$q_1^-$};
\node at (2.5,-0.3){$p_2^-$};
\node at (4.3,-0.3){$q_3^-$};
\node at (6.3,-0.3){$p_4^-$};
\node at (8.5,-0.3){$q_5^-$};

      \filldraw[black] (-1,0) circle (2pt);
      \end{tikzpicture}
    \caption{\label{fig:51} The pentafoil knot $5_1$ with strands and crossings labeled }
\end{figure}
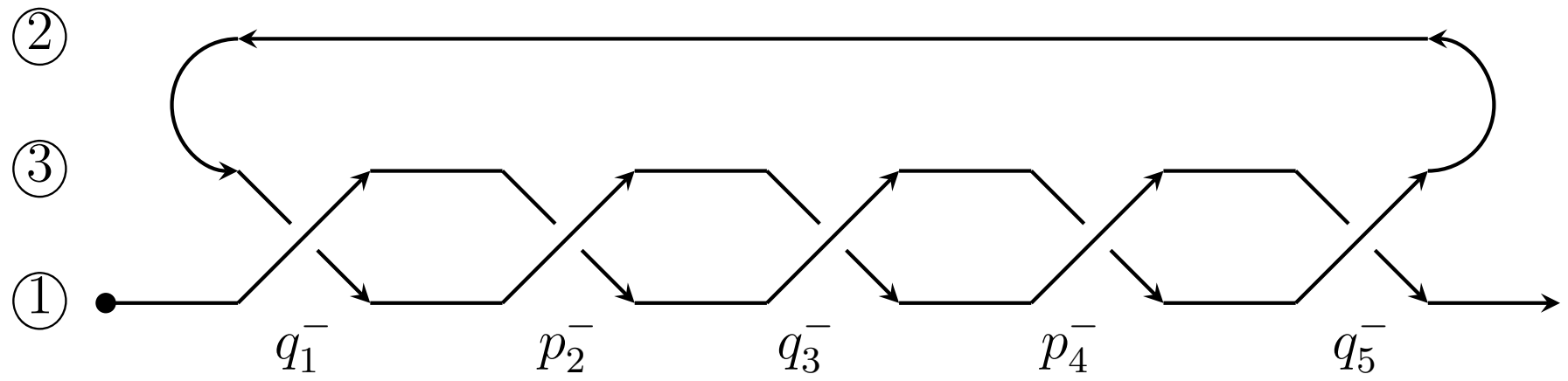

We conclude by summarizing these computations for the base cases in the inductive proof below: the sizes of the partially doubly-counted sets are $t(3)=t(4)=1$ and $t(5)=3$; the total Casson invariants from $T(c)$ are $v(3)=1$, $v(4)=-1$, and $v(5)=7$;  and the first-strand signed ascending sums are $f(3)=-1$, $f(4)=-1$, and $f(5)=-5$.

\subsection{The recursion}
\label{subsec:recursion}  
\begin{figure}[h!]

		\begin{tikzpicture}[scale=0.60]	
				
		\draw[very thick, stealth-] (7,0) arc (-90:90:.5);		
		\draw[very thick, stealth-] (7,1) -- (6,2);
		\draw[very thick, stealth-] (4,2) -- (5.4,0.6);
		
		\draw[very 	thick, -stealth] (5,0) -- (6,1);
		\draw[very thick, -stealth] (4,0) -- (5,0);
		
		\draw[very thick, -stealth] (4,1) -- (4.4,1.4);
		\draw[very thick, stealth-] (5,2) -- (4.6,1.6);
		\draw[very thick, -stealth] (5,2) -- (6,2);
		\draw[very thick] (6,0) -- (7,0);
		\draw[very thick] (5,1) -- (5.4,0.6);
		\draw[very 	thick, stealth-] (5.6,0.4) -- (6,0);
		\node at (4.6,2.6){$p_{c-2}^-$};
		\node at (6.4,2.6){$q_c^+$};
		\node at (5.5,-0.4){};
		\draw[very thick] (6,1) -- (6.4,1.4);
		\draw[very thick, stealth-] (7,2) -- (6.6,1.6);
		\node at (3,1){\textcircled{1}};
		\node at (3,0){\textcircled{3}};
		\node at (8.5,1){$\leadsto$};	
            			\node[above] at (8.5,1){IE};
\begin{scope}[shift={(12,0)}]
        \draw[very thick, -stealth] (1,1) arc (-90:90:.5);
		\node at (-2,1){\textcircled{1}};
		\node at (-2,0){\textcircled{3}};
		\draw[very thick, stealth-] (-1,2) -- (1,2);
		
		\draw[very 	thick, -stealth] (0,0) -- (1,1);
		\draw[very thick, -stealth] (-1,0) -- (0,1);
		\draw[very thick, -stealth] (-1,1) -- (-0.6,0.6);
		\draw[very thick] (0,0) -- (-0.4,0.4);
		\draw[very thick] (0,1) -- (0.4,0.6);
		\draw[very 	thick, -stealth] (0.6,0.4) -- (1,0);
		\node at (-0.65,-0.5){$p^-_{c-2}$};
		\node at (0.65,-0.5){$q^-_{c-1}$};
\end{scope}
\begin{scope}[shift={(15,0)}]
		
		\draw[very thick, -stealth] (4,2) -- (5,2) -- (6,1) -- (6.4,.6);
		
		\draw[very 	thick, -stealth] (4.6,.4) -- (5,0) -- (6,0) -- (7,1);
		\draw[very thick, -stealth] (4,1) -- (4.4,0.6);
		\draw[very thick, stealth-] (4,0) -- (5,1) -- (5.4,1.4);
		\draw[very thick, stealth-] (5.6,1.6) -- (6,2) -- (7,2);
		
		\draw[very thick, -stealth] (7,1) arc (-90:90:.5);
		\draw[very thick, -stealth] (6.6,.4) -- (7,0);

		
		\node at (4.6,-0.5){$p_{c-2}^+$};
		\node at (6.4,-0.5){$q_c^-$};
		\node at (5,1.5){};
		
		\node at (3,2){\textcircled{3}};
		\node at (3,1){\textcircled{1}};
		\node at (8.5,1){$\leadsto$};
            			\node[above] at (8.5,1){IO};
\end{scope}
\begin{scope}[shift={(27,0)}]	
		\node at (-2,2){\textcircled{3}};
		\node at (-2,1){\textcircled{1}};
		\draw[very 	thick, -stealth] (-1,1) -- (-.6,1.4);
		\draw[very 	thick, -stealth] (-.4,1.6) -- (0,2) -- (1,1);

        \draw[very thick, -stealth] (-1,2) -- (0,1) -- (.4,1.4);
		\draw[very thick, -stealth] (.6,1.6) -- (1,2);
		\draw[very thick] (0,1) -- (-0.4,1.4);
		\draw[very thick] (0,2) -- (0.4,1.6);
		\draw[very 	thick, -stealth] (0.6,1.4) -- (1,1);
		
		\draw[very 	thick, -stealth] (1,0) -- (-1,0);
		
		\node at (-0.7,2.5){$p_{c-2}^+$};
		\node at (0.7,2.5){$q_{c-1}^+$};
		\draw[very thick, stealth-] (1,0) arc (-90:90:.5);
\end{scope}
        \end{tikzpicture}

\vspace{.2cm}

		\begin{tikzpicture}[scale=0.60]
\begin{scope}[shift={(-8,0)}]		
		\node at (8,2){\textcircled{3}};
		\node at (8,0){\textcircled{1}};
		
		\draw[very thick, -stealth] (9,0) -- (10.4,1.4);
		\draw[very thick, stealth-] (9,1) -- (9.4,0.6);
		\draw[very thick, stealth-] (9,1) -- (9.4,0.6);
		\draw[very thick] (10,0) -- (9.6,0.4);
		\draw[very thick, -stealth] (9,2) -- (10,2);
		\draw[very thick, stealth-] (11,1) -- (10,2);
		\draw[very thick, stealth-] (11,2) -- (10.6,1.6);
		\draw[very thick, stealth-] (10,0) -- (12,0);
		\draw[very thick, -stealth] (11,2) -- (12,1);
		
		\draw[very thick, -stealth] (11,1) -- (11.4,1.4);
		\draw[very thick, stealth-] (12,2) -- (11.6,1.6);
		
		\node at (10.2,2.55){$p_{c-1}^+$};
		\node at (11.6,2.55){$q_{c}^+$};
		
		\draw[very thick, stealth-] (12,0) arc (-90:90:.5);			
		\node at (13.2,1){$\leadsto$};	
            			\node[above] at (13.2,1){IIE};
\end{scope}
\begin{scope}[shift={(4,0)}]
        \draw[very thick, -stealth] (4,2) -- (5.4,0.6);
		
		\draw[very 	thick, -stealth] (5,0) -- (6,1);
		\draw[very thick, -stealth] (4,0) -- (5,0);
		
		\draw[very thick, stealth-] (4,1) -- (4.4,1.4);
		\draw[very thick, -stealth] (5,2) -- (4.6,1.6);
		\draw[very thick, stealth-] (5,2) -- (6,2);
		
		\draw[very thick] (5,1) -- (5.4,0.6);
		\draw[very 	thick, -stealth] (5.6,0.4) -- (6,0);
		\node at (4.6,2.4){};
		
		\node at (5.5,-0.4){$q_{c-1}^-$};
		\draw[very thick, -stealth] (6,1) arc (-90:90:.5);	
		
		\node at (3,2){\textcircled{3}};
		\node at (3,0){\textcircled{1}};
\end{scope}
\begin{scope}[shift={(12,0)}]

		\draw[very thick, -stealth] (7,1) arc (-90:90:.5);
		\draw[very thick, -stealth] (6.6,0.4) -- (7,0);
		
		\draw[very thick, -stealth] (4,2) -- (5.4,0.6);

		\draw[very 	thick, -stealth] (5,0) -- (6,1);
		\draw[very thick, -stealth] (4,0) -- (5,0);
		\draw[very thick, -stealth] (6,0) -- (7,1);
		\draw[very thick, stealth-] (4,1) -- (4.4,1.4);
		\draw[very thick, -stealth] (5,2) -- (4.6,1.6);
		\draw[very thick, stealth-] (5,2) -- (7,2);
		
		\draw[very thick] (5,1) -- (5.4,0.6);
		\draw[very 	thick, -stealth] (5.6,0.4) -- (6,0);
		\node at (4.6,2.4){};
		\node at (6.6,2.4){};
		\node at (5.4,-0.4){$p_{c-1}^-$};
		\node at (6.6,-0.4){$q_c^-$};
		\draw[very thick] (6,1) -- (6.4,0.6);
		
		\node at (3,2){\textcircled{1}};
		\node at (3,0){\textcircled{3}};

		\node at (8.5,1){$\leadsto$};
           			\node[above] at (8.5,1){IIO};
\end{scope}                    
\begin{scope}[shift={(24,0)}]		
		
		\node at (-2,2){\textcircled{1}};
		\node at (-2,0){\textcircled{3}};
		\draw[very thick, stealth-] (1,1) -- (0,2);
		\draw[very thick, -stealth] (-1,2) -- (0,2);
		\draw[very thick, stealth-] (1,0) arc (-90:90:.5);
		
		\draw[very thick, -stealth] (-1,0) -- (0.4,1.4);
		\draw[very thick, stealth-] (1,2) -- (0.6,1.6);
		\draw[very thick, stealth-] (-1,1) -- (-0.6,0.6);
		\draw[very thick] (0,0) -- (-0.4,0.4);
		\draw[very thick] (0,0) -- (1,0);
		
		\node at (-0.4,-0.5){};
		\node at (0.4,2.55){$q_{c-1}^+$};
\end{scope}		
		\end{tikzpicture}

\vspace{.2cm}

		\begin{tikzpicture}[scale=0.60]
		\begin{scope}[shift={(-7.75,0)}]

		\node at (7,2){\textcircled{1}};
		\node at (7,1){\textcircled{3}};

		\draw[very thick, stealth-] (8,0) -- (11,0);
		
		\draw[very thick, -stealth] (8,2) -- (9,1);
		\draw[very thick, -stealth] (8,1) -- (8.4,1.4);
		\draw[very thick, stealth-] (9,2) -- (8.6,1.6);

		\draw[very thick, -stealth] (9,2) -- (10,1);
		\draw[very thick, -stealth] (9,1) -- (9.4,1.4);
		\draw[very thick, stealth-] (10,2) -- (9.6,1.6);
		
		\draw[very thick, -stealth] (10,2) -- (11,1);
		\draw[very thick, -stealth] (10,1) -- (10.4,1.4);
		\draw[very thick, stealth-] (11,2) -- (10.6,1.6);
		
		\draw[very thick, stealth-] (11,0) arc (-90:90:.5);

		\node at (8.1,2.5){\small $q_{c-2}^+$};
		\node at (9.5,2.5){\small $p_{c-1}^+$};
		\node at (10.7,2.5){\small $q_{c}^+$};
		\node at (12.5,1){$\leadsto$};
           			\node[above] at (12.5,1){IIIE};
\end{scope}
\begin{scope}[shift={(3.25,0)}]
        \draw[very thick, stealth-] (5,0) arc (-90:90:.5);
		
		\draw[very thick, -stealth] (4,2) -- (5,1);

		\draw[very thick, stealth-] (4,0) -- (5,0);
		
		\draw[very thick, -stealth] (4,1) -- (4.4,1.4);
		\draw[very thick, stealth-] (5,2) -- (4.6,1.6);
		
		\node at (3,2){\textcircled{1}};
		\node at (3,1){\textcircled{3}};

		\node at (4.5,2.5){\small $q_{c-2}^+$};
\end{scope}
\begin{scope}[shift={(6.25,0)}]		
		\node at (8,1){\textcircled{3}};
		\node at (8,0){\textcircled{1}};
		
		\draw[very thick, -stealth] (9,0) -- (10,1);
		\draw[very thick, -stealth] (9,1) -- (9.4,0.6);
		\draw[very thick] (10,0) -- (9.6,0.4);
		\draw[very thick, -stealth] (11,1) -- (11.4,0.6);
		\draw[very thick, stealth-] (12,0) -- (11.6,0.4);	
		\draw[very thick, -stealth] (10,0) -- (11,1);
		\draw[very thick, -stealth] (10,1) -- (10.4,0.6);
		\draw[very thick] (11,0) -- (10.6,0.4);
		\draw[very thick] (11,0) -- (12,1);
		\draw[very thick, stealth-] (9,2) -- (12,2);
		\draw[very thick, -stealth] (12,1) arc (-90:90:.5);
		
		\node at (13.5,1){$\leadsto$};
           			\node[above] at (13.5,1){IIIO};
		
		\node at (9.1,-0.5){\small $q_{c-2}^-$};
		\node at (10.5,-0.5){\small $p_{c-1}^-$};
		\node at (11.7,-0.5){\small $q_{c}^-$};
\end{scope}
\begin{scope}[shift={(17.25,0)}]
        \draw[white] (5,0)--(7.25,0);
		
		\draw[very thick, -stealth] (5,1) -- (5.4,0.6);
		
		\draw[very 	thick, -stealth] (5,0) -- (6,1);

		\draw[very thick, stealth-] (5,2) -- (6,2);
		
		\draw[very thick] (5,1) -- (5.4,0.6);
		\draw[very 	thick, -stealth] (5.6,0.4) -- (6,0);
		
		\draw[very thick, -stealth] (6,1) arc (-90:90:.5);
		
		\node at (5.5,-0.5){\small $q_{c-2}^-$};	
		\node at (4,1){\textcircled{3}};
		\node at (4,0){\textcircled{1}};
    \end{scope}
        \end{tikzpicture}
\vspace{.2cm}


		\begin{tikzpicture}[scale=0.60]
\begin{scope}[shift={(-8,0)}]
        \node at (8,2){\textcircled{1}};
		\node at (8,1){\textcircled{3}};
		\draw[very thick, stealth-] (12,0) arc (-90:90:.5);	
		\draw[very thick, stealth-] (9,0) -- (10,1);
		\draw[very thick, -stealth] (9,1) -- (9.4,0.6);
		\draw[very thick] (10,0) -- (9.6,0.4);
		
		\draw[very thick, -stealth] (10,0) -- (11.4,1.4);
		\draw[very thick, stealth-] (10,1) -- (10.4,0.6);
		\draw[very thick] (11,0) -- (10.6,0.4);
		
		\draw[very thick, -stealth] (9,2) -- (11.2,2);
		\draw[very thick, stealth-] (12,1) -- (11,2);
		\draw[very thick, -stealth] (12,0) -- (11,0);
		\draw[very thick, stealth-] (12,2) -- (11.6,1.6);

		\node at (9.5,-0.5){\small $q_{c-2}^+$};
		\node at (11.5,2.5){\small $q_{c}^+$};
		\node at (13.5,1){$\leadsto$};
           			\node[above] at (13.5,1){IVE};
\end{scope}
\begin{scope}[shift={(4,0)}]
		\draw[very thick, stealth-] (5,0) arc (-90:90:.5);	
		
		\draw[very thick, -stealth] (4,2) -- (5,1);
		\node at (3,2){\textcircled{1}};
		\node at (3,1){\textcircled{3}};

		\draw[very thick, stealth-] (4,0) -- (5,0);
		
		\draw[very thick, -stealth] (4,1) -- (4.4,1.4);
		\draw[very thick, stealth-] (5,2) -- (4.6,1.6);
		\node at (4.5,2.5){\small $q_{c-2}^+$};
\end{scope}
\begin{scope}[shift={(7,0)}]
		\node at (8,1){\textcircled{3}};
		\node at (8,0){\textcircled{1}};

		\draw[very thick, -stealth] (9,0) -- (11,0);
		
		\draw[very thick, stealth-] (9,2) -- (10,1);
		\draw[very thick, -stealth] (9,1) -- (9.4,1.4);
		\draw[very thick, stealth-] (10,2) -- (9.6,1.6);
		
		\draw[very thick, stealth-] (12,0) -- (11.6,0.4);
		\draw[very thick, -stealth] (10,2) -- (11.4,0.6);
		\draw[very thick, stealth-] (10,1) -- (10.4,1.4);
		\draw[very thick, -stealth] (11,2) -- (10.6,1.6);
		
		\draw[very thick, stealth-] (11,2) -- (12,2);
		\draw[very thick, stealth-] (12,1) -- (11,0);

		\draw[very thick, -stealth] (12,1) arc (-90:90:.5);		
		\node at (13.5,1){$\leadsto$};		
           			\node[above] at (13.5,1){IVO};
		
		\node at (11.5,-0.5){\small $q_{c}^-$};
		\node at (9.5,2.5){\small $q_{c-2}^-$};

\end{scope}		

        \begin{scope}[shift={(19,0)}]
        \draw[white] (5,0)--(6.5,0);

		\draw[very thick, -stealth] (4,1) -- (4.4,0.6);
		
		\draw[very 	thick, -stealth] (4,0) -- (5,1);

		\draw[very thick, stealth-] (4,2) -- (5,2);
		\node at (4.5,-0.5){\small $q_{c-2}^-$};
		\draw[very thick] (4,1) -- (4.4,0.6);
		\draw[very 	thick, -stealth] (4.6,0.4) -- (5,0);
		
		\draw[very thick, -stealth] (5,1) arc (-90:90:.5);	
		\node at (3,1){\textcircled{3}};
		\node at (3,0){\textcircled{1}};
        \end{scope}
        \end{tikzpicture}
\caption{\label{fig:cases} The four cases in either parity showcasing the original tail and the replacement tail, where crossings are labeled only when they are ascending on the first strand or descending on the third}
\end{figure}
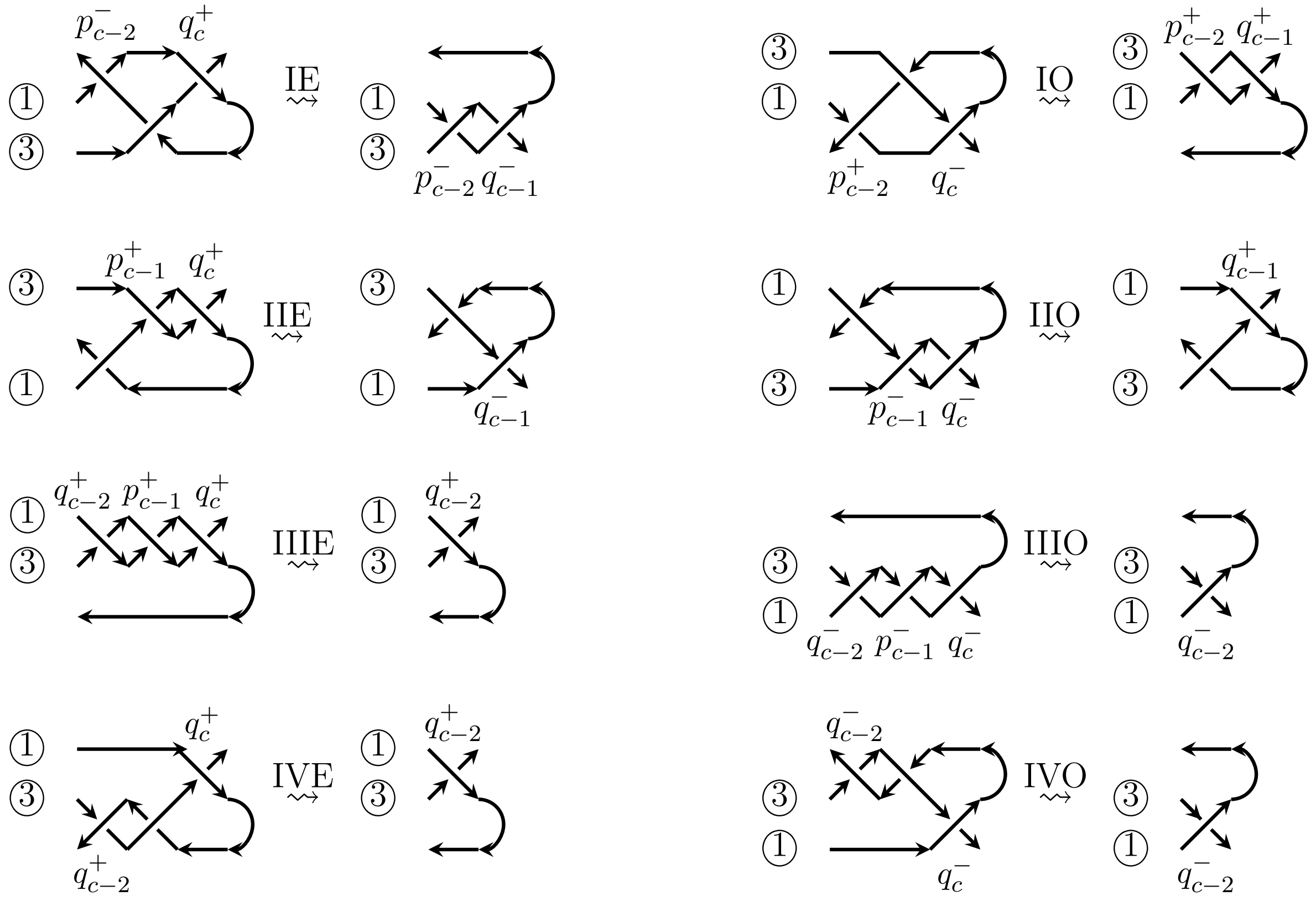

Recall Equations \eqref{eq:tc} and \eqref{eq:tic} from \cite{CohLow} mentioned in Section \ref{sec:background}, 
as well as the four cases of the recursion.  We illustrate the \textbf{tail replacements} for these cases in Figure \ref{fig:cases} and decorate relevant crossings using Notation \ref{not:pq}.

We now use this figure to determine how both the first-strand signed ascending sum and the Casson invariant change under this recursion.

\begin{proposition}
    \label{prop:f+v}
For a word $w\in T(c)$, we consider its image under the bijection described above to be either $w'\in T(c-1)$ or $w''\in T(c-2)$, and we let $D$, $D'$, and $D''$ be the corresponding knot diagrams, respectively.  

The first-strand signed ascending sum satisfies the recursions:
\begin{itemize}
    \item $f(D)=f(D')$ \hspace{20mm} for $w\in T_1(c)$,     
    \item $f(D)=f(D')+(-1)^c$ \hspace{5mm} for $w\in T_2(c)$, 
    \item $f(D)=f(D'')+(-1)^c$ \hspace{5mm} for $w\in T_3(c)$, and 
    \item $f(D)=f(D'')$ \hspace{20mm} for $w\in T_4(c)$.     
\end{itemize}

The Casson invariant satisfies the recursions:
\begin{itemize}
    \item $v_2(D) = v_2(D') + 2(-1)^c f(D')$ \hspace{15mm} for $w\in T_1(c)$,     
   \item $v_2(D) = v_2(D') + 2(-1)^c f(D')+1$ \hspace{8mm} for $w\in T_2(c)$,     
   \item $v_2(D) = v_2(D'') +(-1)^c f(D'')+1$ \hspace{8mm} for $w\in T_3(c)$, and
   \item $v_2(D) = v_2(D'') + (-1)^c f(D'')$ \hspace{15mm} for $w\in T_4(c)$. 
\end{itemize}
\end{proposition}

\begin{proof}
By Lemma \ref{lem:penultimateorder} between the penultimate and last crossings of a knot in $T(c)$, the first strand is at the top for even $c$ and is at the bottom for odd $c$.  We trace these backward through two crossings as shown in Figure \ref{fig:cases}.  Observe that the difference in the first-strand signed ascending sums comes only from the replacement in the tails. 

We start with case I; first assume the crossing number $c$ is even.  Observe that in the tail there is one ascending crossing $p_{c-2}^-$ on strand \textcircled{1} before the replacement and that in the new tail there is still one ascending crossing $p_{c-2}^-$ on strand \textcircled{1} after the replacement; both have a negative sign.  In case I when $c$ is odd we count $p_{c-1}^+$ before and $p_{c-1}^+$ after, both with the same sign.  Thus $f(D)=f(D')$ for $w\in T_1(c)$ for either parity.

In both cases II and III for $c$ even we count $p_{c-1}^+$ before and nothing in the new tail.  For $c$ odd we count $p_{c-1}^-$ before and nothing in the new tail.  Thus for either parity we obtain $f(D)=f(D')+(-1)^c$ for $w\in T_2(c)$ and $f(D)=f(D'')+(-1)^c$ for $w\in T_3(c)$.

In case IV for either parity no ascending crossing appears in either the original or new tails, and so $f(D)=f(D'')$ for $w\in T_4(c)$.

Next we consider the total Casson invariant in these four cases by investigating valid pairs as in Proposition \ref{prop:cassonpairs}.  
By Lemma \ref{lem:validcrossings}  we only need consider ascending crossings on the first strand \textcircled{1} and descending crossings on the third \textcircled{3}.  

In case I, both the original tail and the replacement tail have a single relevant ascending crossing and a single relevant descending crossing.  Thus every valid pair contributing to $v_2(D)$ will also contribute to $v_2(D')$ and vice versa.  The only difference is that for even $c$ the last descending crossing $q_{c}^+$ is replaced by a descending crossing $q_{c-1}^-$ with opposite sign and that for odd $c$ the last descending crossing $q_{c}^-$ is replaced by a descending crossing $q_{c-1}^+$ with opposite sign.  Thus the valid pairs containing this last descending crossing will contribute the opposite sign after the replacement.  For even $c$ these pairs contribute $+f(D)=f(D')$ originally and $-f(D)=-f(D')$ after the replacement, motivating the introduction of this count in Defintion \ref{def:f}, and for odd $c$ it is vice versa.  This yields $v_2(D)=v_2(D')+2(-1)^cf(D')$ for $w\in T_1(c)$.

In case II we have a similar argument since for even $c$ the last descending crossing $q_c^+$ is replaced with $q_{c-1}^-$ of opposite sign and for odd $c$ the last $q_c^-$ is replaced with $q_{c-1}^+$ of opposite sign.  Additionally there is a single valid pair $(p_{c-1}^+,q_c^+)$ for even $c$ or $(p_{c-1}^-,q_c^-)$ for odd $c$ that does not appear in the replacement tail.  Together these yield $v_2(D)=v_2(D')+2(-1)^cf(D')+1$ for $w\in T_2(c)$.

In case III there is another single valid pair $(p_{c-1}^+,q_c^+)$ for even $c$ or $(p_{c-1}^-,q_c^-)$ for odd $c$ that does not appear in the replacement tail.  There is also an additional descending crossing $q_{c-2}^+$ for even $c$ or $q_{c-2}^-$ for odd $c$ in the original tail so that we get an additional $(-1)^cf(D'')$ contributing to $v_2(D)$ but not to $v_2(D'')$. Together these yield $v_2(D)=v_2(D'')+(-1)^cf(D'')+1$ for $w\in T_3(c)$.

In case IV we have a similar argument for an additional descending crossing $q_{c-2}^+$ for even $c$ or $q_{c-2}^-$ for odd $c$ in the original tail, yielding $v_2(D)=v_2(D'')+(-1)^cf(D'')$ for $w\in T_4(c)$.
\end{proof}

Since these four cases $T_i(c)$ partition $T(c)$, we can add together the contributions from Proposition \ref{prop:f+v} to obtain the following recursions.

\begin{lemma}
\label{lem:f+v}
The first-strand signed ascending sum for crossing number $c$ from Definition \ref{def:f} satisfies the recursion
\begin{equation}
\label{eq:frec}
    f(c)=f(c-1)+2f(c-2)+2(-1)^c t(c-2),
\end{equation}
and the total Casson invariant from $T(c)$ from the introduction   
satisfies the recursion
\begin{equation}
\label{eq:vrec}
    v(c)= v(c-1) + 2v(c-2) + 2(-1)^c[f(c-1)+f(c-2)]+2t(c-2)
\end{equation}
\end{lemma}

\begin{proof}
To prove Equation \eqref{eq:frec}, first recall that cases I and II partition $T(c-1)$ so that summing $f(D')$ for all $w\in T_1(c)$ and $f(D')$ for all $w\in T_2(c)$ total to $f(c-1)$.  Cases III and IV together give two copies of $f(c-2)$ summing over all $w\in T_3(c) \sqcup T_4(c)$.  This leaves the contribution of $(-1)^c$ for each $w\in T_2(c)\sqcup T_3(c)$, and by Equation \eqref{eq:tic}, we have $|T_2(c)|=t(c-2)=|T_3(c)|$.

To prove Equation \eqref{eq:vrec}, we use the same argument as above to first get $v(c-1)+2v(c-2)$, and we apply it again to obtain $2(-1)^c f(c-1)+2(-1)^c f(c-2)$.  This leaves the contribution of $+1$ for each $w\in T_2(c)\sqcup T_3(c)$, which additionally gives $2t(c-2)$.
\end{proof}

Applying these recursions to our initial values summarized at the end of Subsection \ref{subsec:examples} gives us the first-strand signed ascending sum $f(6)=f(5)+2f(4)+2t(4)=-5-2+2=-5$ and total Casson invariant $v(6)=v(5)+2v(4)+2[f(5)+f(4)]+2t(4)=7-2+2[-5-1]+2=-5$.

\subsection{Computations}
\label{subsec:computations}
Equation \eqref{eq:vrec} and some foresight below in Lemmas \ref{lem:computations} and  \ref{lem:w}  
lead us to define a few auxiliary functions:
\begin{align*}
g(c) \; & =\;  (-1)^c[f(c)+f(c-1)] \text{ for $c\geq 4$,} \\
h(c) \; & =\; g(c)-2t(c-3)\text{ for $c\geq6$, and }\\
w(c) \; & =\; v(c) + v(c-1) \text{ for $c\geq 4$.}
\end{align*}

 Note that the summary of small values of $f$ and $v$ at the end of Subsection \ref{subsec:examples} and the recursion at the end of Subsection \ref{subsec:recursion} give us the initial value $g(4)=+[f(4)+f(3)]=(-1+-1)=-2$, the value $g(5)=-[f(5)+f(4)]=-(-5+-1)=6$, and finally the 
value $g(6)=+[f(6)+f(5)]=-10$.  This in turn leads us to the initial value of $h(6)=g(6)-2t(3)=-10-2=-12$.  We also have $w(4)=0$, $w(5)=6$, and $w(6)=2$.

\begin{lemma}
	\label{lem:computations}
	The term $g(c)=(-1)^c[f(c)+f(c-1)]$ can be presented in closed form by:
	\begin{equation}\label{eq:g}
	g(c)       \; =\; \frac{1}{3}2^{c-4} -3(-2)^{c-4} + \frac{2}{3}(-1)^{c}.    
	\end{equation}
\end{lemma}
\begin{proof}  Adding $f(c-1)$ to both sides of Equation \eqref{eq:frec} in Lemma \ref{lem:f+v}, we get
	\begin{align}
	\notag f(c)+f(c-1) \; & =\; 2[f(c-1)+f(c-2)]+2(-1)^ct(c-2) \\
	\notag \underbrace{(-1)^c[f(c)+f(c-1)]} \; & =\; 2(-1)^c[f(c-1)+f(c-2)]+2t(c-2) \\
	g(c) \hspace{14mm}   \; & =\; -2g(c-1)+2t(c-2) \label{eq:2} \\
	\notag     \; & =\; -2g(c-1)+2[t(c-3)+2t(c-4)]\ \ \ \text{using Equation \eqref{eq:tc}}\\
	\notag \underbrace{g(c) - 2t(c-3)}   \; & =\; -2g(c-1)+4t(c-4) \\
	h(c) \hspace{9mm}        \; & =\; -2 h(c-1) \label{eq:3}\\
	\notag        \; & =\; (-2)^{c-6}h(6) \text{ by Equation \eqref{eq:3} repeatedly } \\
	\notag        \; & =\; (-2)^{c-6}(-12) = -3(-2)^{c-4}.
	\end{align}
Thus $g(c)= h(c)+2t(c-3) 
	=[-3(-2)^{c-4}] + [\frac{2}{3}2^{c-5}-\frac{2}{3}(-1)^{c-3}]$ using Equation \eqref{eq:tc}.
\end{proof}

\begin{remark}
	\label{rem:g}
	Equation \eqref{eq:g} gives for $g(c)$ the sequence  $-2,6,-10,26,-42,106,-170,\ldots$ starting at $c=4$. We observe that $|g(c)|$ aligns with A061547 from the Online Encyclopedia of Integer Sequences  \cite{OEIS61547},  which counts the number of derangements of $\{1,\ldots,n\}$ avoiding $132$ and $213$.
\end{remark}

\begin{lemma}\label{lem:w} For $c\geq 4$  the auxiliary function $w(c)=v(c)+v(c-1)$ has the closed form
	\begin{equation}
	\label{eq:w}
	w(c) \;  = \; \frac{1}{3}c2^{c-4}  -\frac{1}{9}2^{c-5} -3(-1)^{c}2^{c-5} + \frac{2}{9}(-1)^c 
	\end{equation}
	or equivalently
	\begin{equation*}
	\label{eq:w2}
	w(c) \;  = \; \frac{1}{48}c2^{c}  - \frac{1}{288}2^{c} -\frac{3}{32}(-1)^{c}2^c + \frac{2}{9}(-1)^c. 
	\end{equation*}
\end{lemma}

\begin{proof}
	We add $v(c-1)$ to both hand sides of Equation \eqref{eq:vrec} in Lemma \ref{lem:f+v} and derive 
	\begin{align*}
	\notag 	w(c) \; & = \; 2w(c-1)+[-2g(c-1)+2t(c-2)]\\
	\notag 		\; & = \; 2w(c-1)+g(c)\ \text{by Equation \eqref{eq:2}}.
	\end{align*}

We apply this recursion repeatedly and use Equation \eqref{eq:g} and calculate to obtain 
        \begin{align*}
	w(c) \; & = \; 2^{c-5}w(4)+\sum_{i=0}^{c-5}2^i g(c-i)\\ 
    \; & = \; 2^{c-5}\cdot 0+\sum_{i=0}^{c-5}2^i \left[ \frac{1}{3}2^{c-i-4}-3(-2)^{c-i-4}+\frac{2}{3}(-1)^{c-i}\right]\\
    \; & = \; \frac{1}{3}\sum_{i=0}^{c-5}2^i 2^{c-i-4} -3 \sum_{i=0}^{c-5}2^i(-2)^{c-i-4}+\frac{2}{3}\sum_{i=0}^{c-5}2^i(-1)^{c-i}\\
    \; & = \; \frac{1}{3}\sum_{i=0}^{c-5}2^{c-4} -3 \sum_{i=0}^{c-5}(-1)^{c-i-4}2^{c-4}+\frac{2}{3}(-1)^c\sum_{i=0}^{c-5}2^i(-1)^{i}\\
    \; & = \; \frac{1}{3}2^{c-4}(c-4) -3 \cdot 2^{c-4}(-1)^c\sum_{i=0}^{c-5}(-1)^{i}+\frac{2}{3}(-1)^c\sum_{i=0}^{c-5}(-2)^{i}\\
    \; & = \; \frac{1}{3}2^{c-4}(c-4) -3 \cdot 2^{c-4}(-1)^c\left(\frac{1-(-1)^{c-4}}{1-(-1)}\right)+\frac{2}{3}(-1)^c\left(\frac{1-(-2)^{c-4}}{1-(-2)}\right)\\
    \; & = \; \frac{1}{3}c2^{c-4}-\frac{1}{3}2^{c-2} -3 \cdot 2^{c-5}(-1)^c+3 \cdot2^{c-5}+\frac{2}{9}(-1)^c-\frac{1}{9}2^{c-3}\\
    \; & = \; \frac{1}{3}c2^{c-4}+\left(-\frac{1}{3}2^{c-2} -\frac{1}{9}2^{c-3} +3 \cdot2^{c-5}\right)-3 \cdot 2^{c-5}(-1)^c+\frac{2}{9}(-1)^c,
    \end{align*}
yielding the desired result.
\end{proof}

Although we originally derived the formula for $v(c)$ by writing $v(c)=\displaystyle\sum_{i=0}^{c-7}(-1)^iw(c-i)$
and expanding the telescoping series, we instead prove 
Theorem \ref{thm:v} by induction on $c\geq 6$.

\begin{proof}
The base case $v(6)=-5$ from the end of Subsection \ref{subsec:recursion} holds in Equation \eqref{eq:mainv}.

Suppose by induction hypothesis that Equation \eqref{eq:mainv} holds for $c-1$:
\begin{align*}
\notag	v(c-1) \; & = \; \dfrac{1}{72}2^{c-1}(c-1)+\dfrac{1}{432}2^{c-1}-\frac{3}{16}(-2)^{c-1} + \frac{2}{9}(-1)^{c-1}(c-1)   +\frac{5}{27}(-1)^{c-1}\\
\notag \; &  =\dfrac{1}{9}c2^{c-4}-\left[\dfrac{2}{9}-\frac{1}{27}\right]2^{c-5}+3(-1)^c2^{c-5}-\dfrac{2}{9}(-1)^{c}c
+\left[-\dfrac{5}{27}+\frac{2}{9}\right](-1)^{c} \\
\notag \; &  =\dfrac{1}{9}c2^{c-4}-\dfrac{5}{27}2^{c-5}+3(-1)^c2^{c-5}-\dfrac{2}{9}(-1)^{c}c
+\dfrac{1}{27}(-1)^{c}.
\end{align*}
Then by Equation \eqref{eq:w} we have
\begin{align*}
\notag	v(c) \; & = \; w(c)-v(c-1) \\
\notag \; &  =\left[ \frac{1}{3}c2^{c-4}  -\frac{1}{9}2^{c-5} -3(-1)^{c}2^{c-5} + \frac{2}{9}(-1)^c\right] \\
\notag \; & \; \hspace{20mm} -\left[\dfrac{1}{9}2^{c-4}c-\dfrac{5}{27}2^{c-5}+3(-1)^c2^{c-5}-\dfrac{2}{9}(-1)^{c}c
+\dfrac{1}{27}(-1)^{c}\right] \\
\notag \; &  = \left[\frac{1}{3}-\frac{1}{9}\right]c2^{c-4}  -\left[\frac{1}{9}-\frac{5}{27}\right]2^{c-5} -\left[3+3\right](-1)^{c}2^{c-5} +\frac{2}{9}(-1)^c c + \left[\frac{2}{9}-\frac{1}{27}\right](-1)^c \\
\notag \; &  = \; \frac{2}{9}c2^{c-4}  +\frac{2}{27}2^{c-5}   -6(-2)^{c-5} + \frac{2}{9}(-1)^c c   +\frac{5}{27}(-1)^{c} \\
\notag \; &  = \; \frac{1}{72}c2^{c}  +\frac{1}{432}2^{c}   -\frac{3}{16}(-2)^{c} + \frac{2}{9}(-1)^c c   +\frac{5}{27}(-1)^{c},
\end{align*}
achieving Equation \eqref{eq:mainv}.
\end{proof}

\section{Contribution from the palindromic words}
\label{sec:palindromes}

In Section \ref{sec:main} we only discussed the contribution to the average Casson invariant for the partially-doubly-counted set $T(c)$, but if we want an exact average, we should additionally consider the contribution from the palindromic words $T_p(c)$.  However, if we are only interested in the asymptotic behavior, we need only show that the contribution from the palindromic words is insignificant as $c$ approaches infinity.  We prove this by imitating an argument from Section 4 of a recent paper \cite{CohLowMadRaa} by the first author on signature.

\begin{proposition}
\label{prop:ignorepalindromes}
    The average $\overline{v_2}(c)$ Casson invariant of $2$-bridge knots with crossing number $c$ satisfies
    \[\overline{v_2}(c)\approx \frac{v(c)}{2|\mathcal{K}(c)|}.\]
\end{proposition}

To prove this we will need to invoke some results about upper and lower bounds on the Casson invariant of any knot.  The following is mentioned by Abe \cite{Abe:Vas} referring to Naosuke Okuda's Master's thesis in Japanese that cannot be accessed electronically.

\begin{lemma}[Okuda \cite{Okuda}]
    \label{lem:bounds}
    The Casson invariant of a knot whose diagram has $n$ crossings satisfies
$$-\frac{n^2}{16} \leq  v_2(K) \leq \frac{n^2}{8}.$$ 
\end{lemma}

Since we cannot access this proof, we instead provide an earlier result with looser bounds by Polyak and Viro.

\begin{lemma}[Polyak-Viro \cite{PolVir:Cas} Theorem 1.E]
     \label{lem:pvboundslooser}
    The Casson invariant of a knot given by a diagram with $n$ crossings satisfies
$$ | v_2(K) | \leq \left\lfloor\frac{n^2}{8}\right\rfloor.$$
\end{lemma}

Now we have the tools to prove Proposition \ref{prop:ignorepalindromes}:

\begin{proof}
    The first  author with Lowrance \cite[Proposition 3.6]{CohLow} proved the number of palindromic words $t_p(c)$ satisfies
    \[t_p(c)=t\left(\left\lfloor\frac{c+3}{2}\right\rfloor\right)=\frac{1}{3}\left(2^{\left\lfloor\frac{c-1}{2}\right\rfloor}-(-1)^{\left\lfloor\frac{c-1}{2}\right\rfloor}\right),\]
   and therefore by Lemma \ref{lem:pvboundslooser} we have that the total palindromic Casson invariant $v_p(c)$ from the introduction 
   satisfies
    \[
    |v_p(c)|\leq \frac{1}{3}\cdot\left\lfloor\frac{c^2}{8}\right\rfloor \left(2^{\left\lfloor\frac{c-1}{2}\right\rfloor}-(-1)^{\left\lfloor\frac{c-1}{2}\right\rfloor}\right).\]
    Thus Ernst and Sumner's Theorem \ref{thm:ernstsumners} implies that $\displaystyle\lim_{c\to\infty}\frac{v_p(c)}{2|\mathcal{K}(c)|}=0$, and then
    \[\overline{v_2}(c) = \frac{v(c)+v_p(c)}{2|\mathcal{K}(c)|}\approx \frac{v(c)}{2|\mathcal{K}(c)|},\]
    as desired.
\end{proof}

\subsection{Examples} \label{subsec:families} As an addendum to the Okuda result above, here are some more details about the maximal and minimal cases via some families of examples.

\begin{proposition}
    [Alvarez-Labastida \cite{AlvLab}]
    \label{prop:torus}
     The Casson invariant for the torus knot $T(2,c)$ for odd $c$ is $v_2 (T (2, c)) = \frac{c^2 - 1}{8}$.
\end{proposition}

One can use our machinery above to reprove this result and also give an idea why it achieves the maximum Casson invariant over 2-bridge knots of crossing number $c$.

\begin{figure}[h!]
    \begin{tikzpicture}
\draw[thick, -stealth] (-1,0) -- (0,0);
\draw[thick, -stealth] (7,0) -- (8,0);
\draw[thick, stealth-] (0,2) -- (7,2);
\draw[thick, -stealth] (0,2) arc (90:270:.5);
\draw[thick, -stealth] (7,1) arc (-90:90:.5);
\draw[white] (1,2.25) -- (2,2.25);

\foreach \x/ \y in {0/0, 1/0, 2/0, 3/0, 4/0, 5/0, 6/0}
    {
    \draw[thick, stealth-] (\x+1,\y) -- (\x,\y+1);
    \draw[color=white, line width=10] (\x,\y) -- (\x+1,\y+1);
    \draw[thick, -stealth] (\x,\y) -- (\x+1,\y+1);
            }
\filldraw[black] (-1,0) circle (2pt);
\end{tikzpicture}
\caption{\label{fig:torus} The torus knot $\sigma_1^7$}
\end{figure}
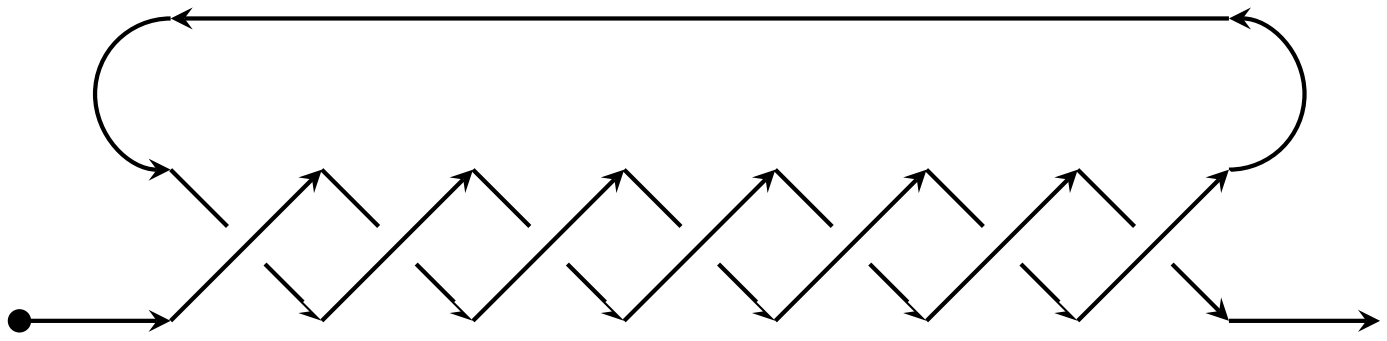

\begin{example}
\label{ex:torus}
Consider the example $T(2,7)$ as in Figure \ref{fig:torus} and any torus knot $T(2,c)$ in general for odd $c$.  Observe that all crossings lie on the first strand as well as the third strand as in Lemma \ref{lem:validcrossings}.  All crossings have the same negative sign, and thus every valid pair will contribute positively.  The odd-indexed crossings are descending, and the even-indexed crossings are ascending.

Requiring that $i(p)<i(q)$ as in our Lemma \ref{lem:validcrossings} above, the contributions to the Casson invariant will come exactly from those pairs with the ascending crossing at position $2j$ and every subsequent descending crossing at positions $2k+1$ for $k\geq j$.  This gives us a sum of 
$$1+2+\ldots+\frac{c-1}{2}=\binom{(c+1)/2}{2}=\frac{1}{2}\cdot\frac{c+1}{2}\cdot\frac{c-1}{2}=\frac{c^2-1}{8}.$$

In this example, we see that all but the first crossing will contribute to the Casson invariant and that all pairs will contribute positively.  See also Subsection \ref{subsec:examples} for the calculations for the smallest torus knots, the trefoil and pentafoil.
\end{example}

\bigskip

\begin{figure}[h!]
\begin{tikzpicture}
\draw[thick, -stealth] (-1,0) -- (0,0);
\draw[thick, -stealth] (6,2) -- (7,2);
\draw[thick, stealth-] (0,2) arc (90:270:.5);
\draw[thick, -stealth] (6,0) arc (-90:90:.5);
\draw[thick, -stealth] (0,2) -- (2,2);
\draw[thick, -stealth] (2,0) -- (6,0);
\foreach \x/ \y in {0/0}
    {
    \draw[thick, -stealth] (\x,\y) -- (\x+1,\y+1);
    \draw[color=white, line width=10] (\x+1,\y) -- (\x,\y+1);
    \draw[thick, -stealth] (\x+1,\y) -- (\x,\y+1);
    }
\foreach \x/ \y in {1/0}
    {
    \draw[thick, stealth-] (\x,\y) -- (\x+1,\y+1);
    \draw[color=white, line width=10] (\x+1,\y) -- (\x,\y+1);
    \draw[thick, stealth-] (\x+1,\y) -- (\x,\y+1);
    }
\foreach \x/ \y in {3/1, 5/1}
    {
    \draw[thick, -stealth] (\x+1,\y) -- (\x,\y+1);
    \draw[color=white, line width=10] (\x,\y) -- (\x+1,\y+1);
    \draw[thick, -stealth] (\x,\y) -- (\x+1,\y+1);
    }
\foreach \x/ \y in {2/1, 4/1}
    {
    \draw[thick, stealth-] (\x+1,\y) -- (\x,\y+1);
    \draw[color=white, line width=10] (\x,\y) -- (\x+1,\y+1);
    \draw[thick, stealth-] (\x,\y) -- (\x+1,\y+1);
    }
\filldraw[black] (-1,0) circle (2pt);
\end{tikzpicture}\\
\begin{tikzpicture}
\draw[thick, -stealth] (-1,0) -- (0,0);
\draw[thick, -stealth] (8,2) -- (9,2);
\draw[thick, stealth-] (0,2) arc (90:270:.5);
\draw[thick, -stealth] (8,0) arc (-90:90:.5);
\draw[thick, -stealth] (0,2) -- (4,2);
\draw[thick, -stealth] (4,0) -- (8,0);
\foreach \x/ \y in {0/0, 2/0}
    {
    \draw[thick, -stealth] (\x,\y) -- (\x+1,\y+1);
    \draw[color=white, line width=10] (\x+1,\y) -- (\x,\y+1);
    \draw[thick, -stealth] (\x+1,\y) -- (\x,\y+1);
    }
\foreach \x/ \y in {1/0, 3/0}
    {
    \draw[thick, stealth-] (\x,\y) -- (\x+1,\y+1);
    \draw[color=white, line width=10] (\x+1,\y) -- (\x,\y+1);
    \draw[thick, stealth-] (\x+1,\y) -- (\x,\y+1);
    }
\foreach \x/ \y in {5/1, 7/1}
    {
    \draw[thick, -stealth] (\x+1,\y) -- (\x,\y+1);
    \draw[color=white, line width=10] (\x,\y) -- (\x+1,\y+1);
    \draw[thick, -stealth] (\x,\y) -- (\x+1,\y+1);
    }
\foreach \x/ \y in {4/1, 6/1}
    {
    \draw[thick, stealth-] (\x+1,\y) -- (\x,\y+1);
    \draw[color=white, line width=10] (\x,\y) -- (\x+1,\y+1);
    \draw[thick, stealth-] (\x,\y) -- (\x+1,\y+1);
    }
\filldraw[black] (-1,0) circle (2pt);
\end{tikzpicture}
\caption{\label{fig:doubletwist} The double twist knot $D(2,4)$, which is $6_1$, written as $\sigma_1^2\sigma_2^{-4}$, and below it the double twist knot $D(4,4)$, which is $8_3$, written as $\sigma_1^4\sigma_2^{-4}$
} 
\end{figure}
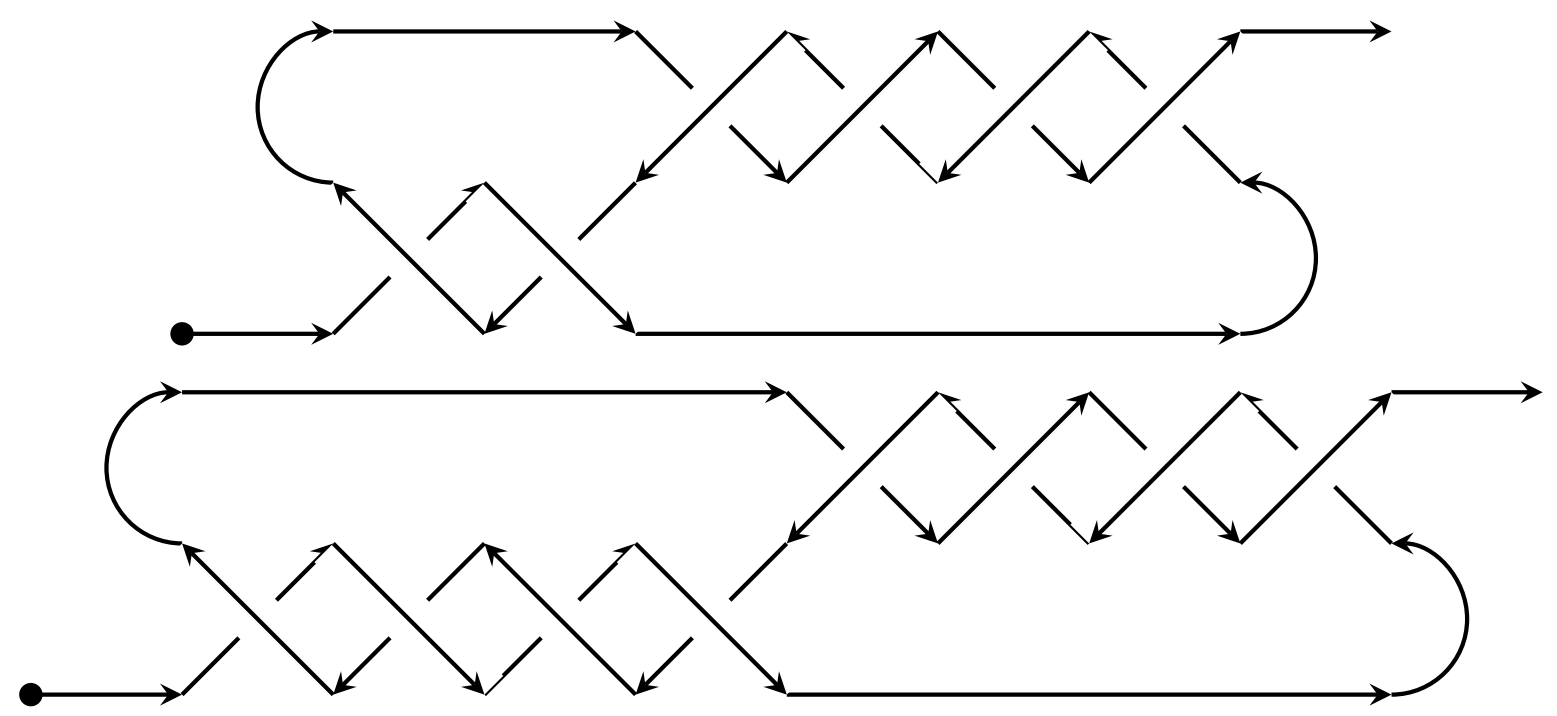

On the other hand, the minimum Casson invariant over our 2-bridge knots is achieved by a different class of knots called the \textbf{double twist knots} $D(a,b)$.  These have just two twist regions giving the knot $\sigma_1^a\sigma_2^{-b}$.  We will restrict our attention to two classes:  $D(2m,2m)$ for $c=4m$ and $D(2m,2m+2)$ for $c=4m+2$ depending on the crossing number modulo 4.  The first such knots are the figure eight knot $4_1$,  the stevedore knot $6_1$, $8_3$, $10_3$, and $12a1287$.  See Figure \ref{fig:doubletwist} for the first two examples.

\begin{proposition}
    \label{prop:doubletwist} 
The Casson invariant for the two classes of double twist knots $D(2m,2m)$ and $D(2m,2m+2)$ are $v_2(D(2m,2m))=-\frac{c^2}{16}$ and $v_2(D(2m,2m+2))=-\frac{c^2-4}{16}$.
\end{proposition}

\begin{proof}
Consider the examples $D(2,4)$ and $D(4,4)$ as in Figure \ref{fig:doubletwist} and any double twist knot $D(2m,2m+2)$ or $D(2m,2m)$ 
in general.  Observe that all crossings are oriented vertically, with both arcs oriented upwards or both downwards; this holds for general $m$.  Observe that all $2m$ crossings in the first twist region lie on the first strand and that all the crossings in the second twist region lie on the third strand as in Lemma \ref{lem:validcrossings}.  All of the crossings in the first twist region have the same positive sign, and all of the crossings in the second twist region have the same negative sign, and thus every pair will contribute negatively.  In the first twist region, the odd-numbered crossings are descending and the even-numbered ones are ascending; in the second twist region, the odd-numbered crossings are ascending and the even-numbered ones are descending.

A valid pair cannot have both crossings from the same twist region because they would not be alternating.  The contributions to the Casson invariant will come exactly from those pairs with any ascending crossing from the first twist region and any descending crossing from the second twist region.  Half of the first $2m$ crossings are ascending, and half of the next $2m$ or $2m+2$ crossings are descending.  This gives us a sum of 
$-m(m)=-\frac{c}{4}\cdot\frac{c}{4}=-\frac{c^2}{16}$ or $-m(m+1)=-\frac{c-2}{4}\cdot\frac{c+2}{4}=-\frac{c^2-4}{16}$ for the two cases.
\end{proof}

In this example, we optimized the number of pairs of crossings that contribute to the Casson invariant where all such pairs contributed negatively.

\begin{remark}
    \label{rem:doubletwist}
Technically these alternating diagrams in Figure \ref{fig:doubletwist} and Proposition \ref{prop:doubletwist} do not result directly from words in $T(c)$ because they have reduced length congruent to 0 modulo 3 as in Definition \ref{def:tc}.  However, these $D(2m,2m)$ can be written as $\sigma_1\sigma_2^{2m-1}\sigma_1^{2m-1}\sigma_2^{-1}$
and these $D(2m,2m+2)$ can be written $\sigma_1\sigma_2^{2m-1}\sigma_1^{2m+1}\sigma_2^{-1}$ coming from words in $T(c)$ that have reduced length congruent to 1 modulo 3.  

As an example, the stevedore knot $6_1$ written $\sigma_1\sigma_2^{-1}\sigma_1^3\sigma_2^{-1}$ shown in Figure \ref{fig:61} is the same as the stevedore knot $6_1$ written $\sigma_1^2\sigma_2^{-4}$ as in Figure \ref{fig:doubletwist}.

This knot equivalence, switching every upper crossing to a lower crossing except for the ones on the end, has been known to the first author since \cite{CoEZKr}.  
\end{remark}

\bibliographystyle{amsalpha}

\newcommand{\etalchar}[1]{$^{#1}$}
\def\cprime{$'$}
\providecommand{\bysame}{\leavevmode\hbox to3em{\hrulefill}\thinspace}
\providecommand{\MR}{\relax\ifhmode\unskip\space\fi MR }
\providecommand{\MRhref}[2]{%
  \href{http://www.ams.org/mathscinet-getitem?mr=#1}{#2}
}
\providecommand{\href}[2]{#2}

\end{document}